\input amstex
\documentstyle{amsppt}

\loadeufb
\loadeusb
\loadeufm
\loadeurb
\loadeusm

\magnification =\magstep 1
\refstyle{A}
\NoRunningHeads

\topmatter
\title  Metrics on universal covering of projective variety 
\endtitle
\author Robert  Treger \endauthor
\address Princeton, NJ 08540  \endaddress
\email roberttreger117{\@}gmail.com \endemail
\endtopmatter

\document

\head
1.  Introduction
\endhead

Let $\phi: X \hookrightarrow  \bold P^{r}$ be a nonsingular connected projective variety of dimension $n\geq 1$. Let $U_X$ denote its universal covering.  Recall that the fundamental group $\pi_1(X)$ is {\it large}\/ if and only if $U_X$ contains no proper holomorphic subsets of positive dimension \cite{Kol}.

Throughout the note we assume that $\pi_1(X)$ {\it large}\/ 
 and {\it residually finite}\/, and the genus $g(C)$ of a general curvilinear section $C\subset X$ is at least 2.
Let $\Cal L=\Cal L_X$ denote  the very ample bundle defining the map $\phi$.

In Section 3, we construct the metric  $\Lambda_\Cal L$ on $U_X$. As an application, in the Appendix, we reproduce a prove of a conjecture of Shafarevich on holomorphic convexity when $\pi_1(X)$ is residually finite (see \cite{T2}).

In Section 4, we assume that $\pi_1(X)$ is, in addition,  nonamenable  (see, e.g., \cite{LS, p.\;300}). We construct a Bergman-type metric with a weight, denoted by $\Sigma_\Cal L$, employing $L^2$ holomorphic functions on $U_X$ 
(the volume form is $dv_\Lambda$).
As a corollary, we obtain that the canonical bundle on $X$, denoted by $\eusm K_X$, is ample.

In Section 5, assuming $\pi_1(X)$ is  nonamenable, we construct another Bergman-type metric with a weight, denoted by $\beta$, 
employing $L^2$ sections of $\eusm K_X^q$ where $q$ is a fixed large positive integer. We obtain a natural embedding into an infinite-dimensional projective space. This means $\pi_1(X)$ is {\it very}\/ large (see \cite{T1}).

In 2010, Campana asked the author whether one can establish the uniformization theorem in \cite{T1} with assumptions on the fundamental group only. 
Since $\pi_1(X)$ is very large provided it is nonamenable, we can employ the proof of our earlier uniformization \cite{T1, Propositions 3.9, 4.2} to obtain the following

\proclaim {Theorem (Uniformization)} Let $X \hookrightarrow  \bold P^{r}$  be a nonsingular connected projective variety of dimension $n$ with {\it large} and {\it residually finite} fundamental group $\pi_1(X)$. If $\pi_1(X)$ is nonamenable then $U_X$ is a bounded Stein domain in $\bold C^n$.
\endproclaim

One of the main ingredients in our proof of the uniformization \cite{T1} was the well-known  Griffiths unformization theorem \cite{G} that says that every algebraic manifold has plenty
of Zarisky open sets which are quotients of bounded domains. 

In Section 5.3, we present another proof of the theorem without appealing to \cite{T1}.

\medskip

In complex geometry, the uniformization problem  is about uniformizing complex analytic functions by means of automorphic functions. 
For a discussion of the problem, see Poincar\'e \cite{Po}, Weierstrass \cite{We, pp.\;95, 232, 304} and Hilbert's 10th problem on the short list \cite{H}.  The problem is to find the assumptions on $X$ such that its universal covering  is a bounded domain. For  modern discussions of the problem, see Siegel \cite {S2, Chap.\;6, Sections 1, 2, pp.\;106-119}, Shokurov \cite{Sho,\;Chap.\;1, $\S$5.2} as well as Piatetski-Shapiro \cite{Pi} and  Shafarevich \cite{PS, Sect.\;4.1}, and an extensive report by
Bers \cite{Be, pp.\;559-609}.

For a review of complex differential geometric approach to the uniformization problem, see a survey by Yau \cite{Y2}.

\bigskip

\head
2.  Preliminaries
\endhead

(2.1) A countable group $\Gamma$ is called {\it amenable}\/ if there is on $\Gamma$ a finite additive, translation invariant nonnegative probability measure (defined for all subsets of $\Gamma$). Otherwise, $\Gamma$ is called {\it nonamenable}.
\medskip

(2.2) {\it Di\'astasis}.\/ The diastasis was introduced by Calabi \cite{C, Chap.\;2}.
Let $M$ denote a complex manifold with a real analytic Kahler metric. Let $\Phi$ denote a real analytic potential of the metric defined in a small neighborhood  $\Cal V \subset M$. Let $z=(z_1,\dots,z_n)$ be a coordinate system in $\Cal V$ and $\bar z=(\bar z_1,\dots,\bar z_n)$ a coordinate system in its conjugate neighborhood $\bar\Cal V\subset \bar M$. 
Let $(\bold p,\bold p)$ be a point on the diagonal of $M \times {\bar M}$ such that the neighborhood $\Cal V \times {\bar \Cal V} \subset M \times {\bar M}$ contains the point.

There exists a unique {\it holomorphic}\/ function $F$ on an open neighborhood of $(\bold p,\bold p)$ such that $F_{(\bold p,\bold p)}= \Phi_\bold p$.
Here $\Phi_\bold p$ is the germ at $\bold p\in M$ of our real analytic function, and $F_{(\bold p,\bold p)}$ is the germ of the corresponding holomorphic function (complexification of $\Phi_\bold p$ \cite {C, Chap.\;2}, \cite{U, Appendix}). 

One considers the sheaf $\Cal A^\bold R_M$ of germs of real analytic functions on $M$, and the sheaf $\Cal A^\bold C_{M \times \bar M}$ of germs of complex holomorphic functions on $M \times \bar M$. For each $\bold p\in M$, we get a natural inclusion $\Cal A^\bold R_{M,\bold p} \hookrightarrow \Cal A^\bold C_{M \times \bar M,(\bold p,\bold p)}$, called a complexification. The above equality is understood in this sense.

  Now, let $p$ and $q$ be two {\it arbitrary}\/ points of $\Cal V$ with coordinates $z(p)$ and $z(q)$. Let $ F(z(p), \overline{z( q)})$ denote the  complex holomorphic function on $\Cal V \times \bar \Cal V$ obtained from from $\Phi$. 
The {\it functional element of diastasis}\/ is defined as follows \cite{C, (5)}:
$$
D_M(p ,q):=F(z(p), \overline{z( p)})\,+\,F(z(q), \overline{z( q)})\,-\,F(z(p), \overline{z( q)})\,-\,F(z(q), \overline{z(p)}). \tag{2.2.1}
$$
We get the germ $D_M(p,q)\in \Cal A^\bold C_{M \times \bar M,(p,q)}$, and $D_M(p,q)$ is uniquely determined  by the Kahler metric, symmetric in $p$ and $q$ and real valued  \cite{C, Prop.\;1,\,2}.
The {\it real}\/ analytic function generated by the above functional element is called the diastasis \cite{C, p.\;3}. 
The diastasis approximates the square of the geodesic distance in the {\it small}\/ \cite{C, p.\;4}.
For $\bold C^r$ with its unitary coordinates,
$
D_{\bold C^r}(p,q) = \sum^r_{i=1} |z_i(q) -z_i(p)|^2.
$

The {\it fundamental property of the diastasis}\/ is that it is inductive on complex submanifolds \cite{C, Chap.\;2, Prop.\;6}.  

Now, let $\bold Q\in M$ be a {\it fixed}\/ point, and $z=(z_1, \dots, z_n)$ a local coordinate system in a small neighborhood $\Cal V_\bold Q \subset M$ with origin at $\bold Q$. The real analytic function $\tilde \Phi_\bold Q (z(p), \overline{z(p)}):= D_M(\bold Q,p)$ on $\Cal V_\bold Q$ is called the {\it  diastasic potential at}\/ $\bold Q$ of the Kahler metric. It is strictly plurisubharmonic function in $p$ \cite{C, Chap.\;2, Prop.\;4}. 

The {\it prolongation}\/ over $M$ of  the germ of   $\tilde \Phi_\bold Q (z(p), \overline{z(p)})$ at $\bold Q$ is a {\it function}\/ $\bold P_M:=\bold P_{M,\bold Q} \in  H^0(\Cal A^\bold R_M, M)$ such that,  for every $u\in M$,  $\bold P_{M}(u)$ coincides  with  $D_M(\bold Q, u) $ meaning $D_M(\bold Q, u) $, initially defined  in a neighborhood of $\bold Q$, can be extended over the {\it whole} $M$. Moreover, the germ of $\bold P_M$ at $u$ is the  diastasic potential of our metric at $u$. 

The  prolongation over $M$ of the germ of  diastasic potential  is not always possible, e.g., there are no strictly plurisubharmonic functions on $\bold P^1$.  
Now, let  $\bold P^N$ be a projective space with the Fubini-Study metric ($1\leq N \leq \infty$). For $\bold Q \in\bold P^N $, we consider Bochner canonical coordinates $z_1, \dots, z_N$ with origin at  $\bold Q$ on the complement of a hyperplane at infinity. By Calabi \cite{Chap.\;4, (27)}, 
$$
	D_{\bold P^r}(\bold Q,p)= \log\bigl(1+\sum^N_{\sigma=1}|z_\sigma(p)|^2 \bigl).
$$
In the homogeneous coordinates $\xi_0,\dots,\xi_N  $, where $z_\sigma := {\xi_\sigma/\xi_0}$, we get
$$
D_{\bold P^N}(\bold Q,p)= \log {\sum^N_{\sigma=0} |\xi_\sigma(p)|^2 \over |\xi_0(p)|^2}.
$$

(2.2.2) {\it Bergman metric.}\/
Let $U\subset \bold C ^n $ be a bounded domain.  
Let $z_1, \dots, z_n$ be a local   coordinate system with origin at a point $\bold Q\in U$. 
Let $\bold B_U$ be the Bergman kernel of $U$.
 By the characteristic property of the diastasic potential (vanishing of some partial derivatives; see  \cite{Bo, pp.\;180-181}, \cite{C, p.\;3,\;p.\;14}, and  \cite{U, Appendix} where this property is explicitly stated):
$$
\partial^{|I|}\bold P_{U, \bold Q}(\bold Q)/\partial z_I=\partial^{|I|}\bold P_{U,\bold Q}(\bold Q)/\partial \bar z_I=0 \; (I:=\{i_1, \dots, i_n\}\; \text{where} \,\; i_1, \dots, i_n \geq 0),
$$
we get  $\log \bold B_U(z,\bar z) - \log \bold B_U(0,0)$ is the diastasic potential at $\bold Q$ of the Bergman metric on $U$. Thus $\bold P_{U,\bold Q}(u(z,\bar z))=\log \bold B_U(z,\bar z)- \log \bold B_U(0,0)$ where $u\in U$, $\bold P_{U,\bold Q}$  is defined over the whole $U$ and $\bold P_{U,\bold Q}(u)= D_U(\bold Q, u)$.

\bigskip

(2.3.1) {\it Tower of coverings}.\/  We consider a
tower of Galois coverings with each $Gal (X_i / X)$ a finite group:
$$
X=X_0 \leftarrow X_1 \leftarrow X_2 \leftarrow \cdots \leftarrow U, \quad \bigcap_i Gal(U/X_i) = \{1\} \; (0\leq i<\infty). \tag{2.3.1.1}
$$
We do not assume $U$ is simply connected. Let $\tau_i$ denote the projection $U\rightarrow X_i$, $\tau:=\tau_0$, and  $\tau_{jk}$ denote the projection $X_j\rightarrow X_k$ $(j\geq k)$.

The hyperplane bundle on $\bold P^{r}$ restricts to $\Cal L_X$,  called a polarization on  $X$ (see, e.g., \cite{Ti, p.\;99}).   Given our polarized Kahler metric  $g$ on  $X$, one can find a Hermitian metric  $h$ on $\Cal L_X$ with its Ricci curvature form equal to the corresponding Kahler form $\omega_g$. 

We consider the volume form of the  Kahler form $\omega_g$. In local coordinates $z_1,\dots,z_n$ on $X$,
$$
dv_g=V_g\! \prod^n_{\alpha=1}\bigl({\sqrt{-1}\over2} \cdot\!dz_\alpha\!\wedge\! d\bar {z}_\alpha \bigl)
$$ 
where $V_g$ is a locally defined positive function. We will employ the same  volume form on all the coverings of $X$. Also, we will employ  the same Hermitian metric on all $\tau^*_{j0}(\Cal L_X)$ and $\tau^*(\Cal L_X)$.
\medskip

(2.3.2) {\it Positive reproducing kernels and Bergman pseudometrics.}\/ The fundamental property of any  Berman-type pseudometric is the existence of a natural {\it continuous} map to a suitable projective space $\bold P(H^*)$ where the corresponding Hilbert space $H$ has a {\it reproducing kernel}.
\smallskip

(2.3.2.1) Let $M$
denote an arbitrary  complex manifold. Let $B(z,w)$ be a
Hermitian positive definite complex-valued  function  on
$M\times M$ which means:
\roster
 \item"{(i)}" $\overline{B(z, w)}=B(w,z), \qquad B(z,z)\geq 0;$
\item"{(ii)}" $\forall \; z_1, \dots, z_N \in M, \quad
 \forall \; a_1, \dots, a_N \in \bold C \implies 
\sum_{j,k}^N B(z_k,z_j)a_j \bar {a}_k \geq 0.$
\endroster
 If $B(z, w)$ is, in addition, holomorphic in the first variable then
$B$ is the reproducing kernel of a {\it unique Hilbert space}\/ $H$ of holomorphic functions on $U$ 
(see Aronszajn \cite{A, p.\;344,(4)} and the articles by Faraut and Kor\'anyi  \cite{FK, pp.\;5-14,
pp.\;187-191}). 
The evaluation at a point $Q\in M,$
$
e_Q:  f\mapsto f(Q),
$
is a continuous linear functional on $H$.
\smallskip

(2.3.2.2) Conversely, given a Hilbert space $H$ ($H\not =0$) of holomorphic functions
on
$M$ with all evaluation maps continuous linear functionals then, by the Riesz
representation theorem, for every $w\in M$ there exists a unique function $B_w\in H$
 such that $f(w)= \langle f, B_w\rangle$ ($\forall f\in H$) and $B(z,w):= B_w(z)$ is the reproducing kernel of
$H$ (which is Hermitian positive definite).

If we assume, {\it in addition}, that $B(z,z) > 0$ for every $z$ then  we can define   
$\log B(z, z)$ and a positive semidefinite  Hermitian form, called the Bergman pseudo-metric
$$
\qquad \qquad ds^2_M =2 \sum g_{j k} dz_j d\overline z_k, \qquad g_{j
k}:={\partial^2\log B(z, z)\over
\partial  z_j\partial\overline z_k  } .
$$
 We get a natural map $
\Upsilon : M \longrightarrow \bold {P}(H^*)$ whose image does not belong to a proper subspace of $\bold {P}(H^*)$   as in \cite{Kob2, Chap.\;4.10, pp.\;224-228}. 

As in  \cite{Kol, Chap.\;7, pp.\;81-84, Lemma-Definition 7.2}, the function $B(z,w)$ can be replaced  by a section of a relevant bundle. 
\smallskip

(2.3.2.3) Let $B_1(z,w)$ and $B_2(z,w)$ be two Hermitian kernels of positive type
(positive matrices or p.\;matrices in the terminology of \cite{A}). We say $B_1 \ll B_2$ if
$B_2(z,w) - B_1(z,w)$ is also a Hermitian kernel of of positive type. 

Recall the following theorem of Aronszajn \cite{A, Part I, Sect.\;7, Theorem II}: If $B$ is the
reproducing kernel of the Hilbert space $F$ with the norm $\parallel\cdot\parallel$, and if
the linear subspace $F_1 \subseteq F$ forms another Hilbert space with the norm 
$\parallel\cdot\parallel_1$, such that  $\parallel f_1\parallel_1 \geq \parallel f_1\parallel$
for every $f_1\in F_1$, then the Hilbert space $F_1$ possesses a reproducing 
kernel $B_1$ satisfying $B_1 \ll B$.

\bigskip

(2.3.3) Assuming $\eusm K_X$ is ample, we fix a large integer $q$ such that for every $i$,  
$\eusm K^q_{X_i}$  is very ample
(see \cite{Kol, 16.5}, \cite{De}). 
The bundle $\eusm K^q_{X_i}$ is equipped with a Hermitian metric $h_{\eusm K^q_{X_i}}:=h^q_{\eusm K_{X_i}}$, where $h_{\eusm K_{X_i}}$ is a Hermitian metric  on $ \eusm K_{X_i}$.
Let $\psi_0,\dots ,\psi_N$ be a  basis of $H^0(X,\eusm K^q)$.
Locally
$\psi_\beta =g_\beta  (z) (dz_1\wedge\cdots \wedge dz_n)^q$. We get an embedding $\sigma: X\hookrightarrow \bold P(H^0(X,K^q)^*)$.
We set 
$$
dv_{X,\eusm K^q}:=\big(\sum_{\beta =0}^N |g_\beta |^2 \big )^{1\over q}(\sqrt{-1})^{n^2} dz_1\wedge\cdots \wedge dz_n \wedge d\bar z_1 \wedge\cdots \wedge d\bar z_n.
$$
This is a volume form  as well as a Hermitian metric on the 
anti-canonical bundle $\eusm K_X^{-1}$.
The associated Ricci form ${\text {Ric}}(dv_{X,\eusm K^q})$ (see, e.g., \cite{Kob2, Chap.\;2.4.4}) is negative. It is known that if we pull back on $X$  the Fubini-Study metric on $\bold P(H^0(X,K^q)^*)$ then its Kahler form differs only by the sign from ${\text {Ric}}(dv_{X,\eusm K^q})$ (see, e.g., \cite{Kob2, Chap.\;7.3}).

\smallskip

(2.3.3.1) Let $\Omega_U^{(n,n)}$  denote the bundle of $(n,n)$-forms. 
As in \cite{Kol, Chap.\;7.1.1.2}, we fix a {\it real}\/ homomorphism 
of bundles:
$$
\Cal H_{U,q} : \eusm K_U^q\otimes \bar \eusm K_U^q \rightarrow \Omega_U^{(n,n)}\simeq\eusm K_U\otimes \bar \eusm K_U. 
$$
If $dv$ is a volume form on $U$ then $h(,\!):={\Cal H_{U,q}/ dv}$ is a Hermitian metric on $\eusm K_U^q $. Given $\Cal H_{U,q}$, we consider the Hilbert space $H=H_{U}$ of  all square-integrable
holomorphic weight $q$ differential forms $\omega$ on $U$. By square-integrable (or $L^2$),  we mean
$$
\int_U\Cal H_{U,q}(\omega\otimes \bar\omega)< \infty.
$$
 
 We assume $H\not=0$. If all evaluation maps are bounded
 then $H$ has a reproducing kernel as in the case of classical Bergman metric \cite{FK, pp.\;8-10, pp.\;187-188}.
Further, if the natural map
$$
 U\longrightarrow \bold P(H^*)
$$
is a holomorphic {\it embedding}\/ then the metric on $U$, induced from $\bold P(H^*)$, is called its $q$-{\it Bergman metric}. It is denoted by $b_{U,q}$ and the corresponding tensor is denoted by $g_{U,\eusm K^q}$. Of course, they depend on the choice of $\Cal H_{U,q}$.

If $W\subsetneq \bold P(H^*)$ is a proper subspace of $\bold P(H^*)$ then the image of $U$ does not lie in $W$ in view of the definition of $H$ in (2.3.3.1), as in the case of classical Bergman metric (\cite{Kob1, Chap.\;7} or \cite{Kob2, Chap.\;4.10, p.\;228}).

Similarly, one defines the Euclidean space $V_i$ of  square-integrable 
holomorphic weight $q$  differential forms on $X_i$.
In the sequel, we assume that each $V_i$ is equipped with a {\it normalized}\/ inner product and a norm
$$
{\parallel\omega \parallel} := \bigg ({1\over Vol(X_i)}{\int_{X_i}\Cal H_{X_i,q}(\omega\otimes \bar\omega)}\bigg )^{\!1/2}\!
$$
(so that we will have the embedding $V_j \hookrightarrow V_k \;(j<k)$ of  Euclidean spaces induced by pullbacks of differential forms provided $\Cal H_{.,q}$ are compatible as in (2.3.3.3) below). 
 Let $\bold P(V^*_i) $ denote the corresponding projective space with its Fubini-Study metric. 
If the natural  map
$$
X_i \longrightarrow \bold P(V^*_i)
$$
is a holomorphic embedding then the induced metric on $X_i$ is called its $q$-{\it Bergman metric}. 
\smallskip

(2.3.3.2) {\it Fundamental domain.}\/ We keep the assumptions of (2.3.1) and (2.3.3)-(2.3.3.1). We set $X_\infty :=U$. We assume each $X_i\;(0\leq i\leq \infty)$ has its $q$-Bergman metric giving an embedding in an appropriate projective space. For $j>i$, let $\Gamma_{ji}:= Gal(X_j/X_i)$
denote the corresponding Galois group.

For $j>i$, a subset $\Cal D \subset X_j$ is called a fundamental domain of $\Gamma_{ji}$ (see, e.g., \cite{Kol, Chap.\; 5.6.2-5.8}) if $X_j =\cup \gamma \Cal D$ and $\gamma \Cal D$ is disjoint from the interior of $\Cal D$ for $\gamma \not = 1$ ($\gamma \in \Gamma_{ji}$).

Now, we pick an arbitrary point $\bold Q_j \in X_j$. We can consider the Dirichlet fundamental domain centered at $\bold Q_j$ (we do not exclude the case $D_{X_j}(p,\bold Q_j)=\infty$ below):
$$
\Cal D_{\bold Q_j}(\Gamma_{ji}):=\big\{ p\in X_j \big|\; D_{X_j}(p,\bold Q_j)\leq D_{X_j}(\gamma(p), \bold Q_j),\; \forall\gamma\in \Gamma_{ji}\}.
$$
Indeed, $ D_{X_j}(q_1,q_2) =  D_{X_j}(\gamma q_1,\gamma q_2)$ for $\forall \gamma\in \Gamma_{ji}$ because of the natural embedding (see (2.3.3)) of ${X_j}$ into the corresponding finite-dimensional or   infinite-dimensional projective space where  $\Gamma_{ji}$ acts by collineations. 

The boundary of $\Cal D_{\bold Q_j}(\Gamma_{ji})$ is a subset of $\Cal D_{\bold Q_j}(\Gamma_{ji})$ where  $\lq\lq\leq"$   is replaced by $\lq\lq="$. The boundary has measure zero (with respect to $dv_g$) because it is at most a countable union of measurable sets of measure zero. For example, the latter can be seen by passing to the complexifications.
\smallskip

(2.3.3.3)  We fix a point $\bold Q \in U$. Let $z=(z_1,\dots,z_n)$ be a coordinate system in a small neighborhood $\Cal V$ with origin at  $\bold Q_i:=\tau_i(\bold Q) \in X_i \; (\forall i)$. 

 Let $\varsigma_q$ be a positive bounded measurable function on $U$.
Also, we assume  $\varsigma_q$ is bounded  away from $0$ on every compact subset of $U$. 

For every $i$, we restrict $\varsigma_q$ on the Dirichlet  fundamental domain of $\Gamma_{\infty i}$ centered at $\bold Q$. We, then, get a bounded positive measurable function on $X_i$. By abuse of notation, we denote the latter function on $X_i$ by the same symbol $\varsigma_q$. 

In the present note, we define a (so-called) {\it compatible  with the tower }\/ (2.3.1.1) sequence of homomorphisms:
$$
\Cal H_{X,q},\; \dots ,\; \Cal H_{X_i,q},\; \dots ,\; \Cal H_{U,q}
$$
as follows.
 In {\it local}\/ coordinates, let $\omega_e=g_e(dz_1\wedge\cdots\wedge dz_n)^q$ be two weight $q$ forms on $U$ or on the Dirichlet fundamental domain of $\Gamma_{\infty i}$ ($e=1, 2$).  Then we define
$$
\Cal H_{.,q}(\omega_1, \bar\omega_2):= (-2\sqrt{-1})^{-n}(-1)^{{n(n-1) \over 2}}\varsigma_q \,g_1 \bar g_2\, dz_1\wedge\cdots \wedge dz_n\wedge d\bar z_1\wedge \cdots \wedge d\bar z_n
$$
provided $\Cal H_{.,q}(\omega_1, \bar\omega_2) $ are well-defined 
(see, e.g., the classical example (2.3.3.4) below). Set $c(n):=(-2\sqrt{-1})^{-n}(-1)^{{n(n-1) \over 2}}$. The homomorphisms are real because 
$$
\sqrt{-1}\cdot2^{-1}dz_\alpha\wedge d\bar z_\alpha = dx_\alpha\wedge dy_\alpha \quad {\text {where}} \quad z_\alpha=x_\alpha+\sqrt{-1}y_\alpha.
$$

We consider $\Cal H_{.,q}$ defined as above with a suitable  $\varsigma_q$.
Recall that we assume all Hilbert spaces are nontrivial. The  corresponding Hilbert spaces have reproducing kernels as in the case of classical weighted Bergman spaces because $\varsigma_q$ is bounded  away from $0$ on every compact subset of $U$  \cite{FK, p.\;10; p.\;188}. 


We pick a point $u\in U$ and consider its image $\tau (u)\in X \subset \bold P^r$. Let $H_\infty \subset \bold P^r$ be a hyperplane at infinity ($\tau (u)\notin H_\infty $). We choose coordinates $z_1, \dots, z_n$ in $X\backslash H_\infty $with origin at $ \tau (u)$. We obtain the same coordinates in a small neighborhood of $u$ in $U$. Let $\xi$ be  a measurable section of $\eusm K_X^q$ such that
$$
\xi|_{X\backslash H_\infty} =(dz_1\wedge \cdots \wedge dz_n)^q.
$$
We denote the inverse image  of $\xi$ on $U$ by the same symbol.

We consider the volume form 
$$
dv_{X,\eusm K^q} =d\mu :=\mu\!\cdot\!(\sqrt{-1})^n dz_1\wedge d\bar z_1 \wedge\cdots \wedge dz_n\wedge d\bar z_n
$$
where $\mu$ is a positive (locally defined) $C^\infty$ function on $X$. Its inverse image on $U$, also denoted by $d\mu$, will be a $\pi_1(X)$-invariant volume form on $U$. Let $\eta_1$ be an arbitrary 
 Hermitian metric on $\eusm K_X$ (\cite{Kol, 5.12, 5.13, 7.1.1} or \cite{Kob2, p.\; 363}). 
Let $\eta_q := \eta_1^q$ be the corresponding Hermitian metric on $\eusm K^q_X$.
By abuse of notation,we denote by $\eta_q$ the  lifting of $\eta_q$ on $\eusm K^q_U$. In our local coordinates, let
$$
\omega = g (dz_1\wedge\cdots \wedge dz_n)^q
$$
be an arbitrary weight $q$ form on $U$.  We will compare 
  $\eta_q (\omega,\bar \omega) d\mu$ and  $\Cal H_{U,q}(\omega\otimes \bar \omega)$.


 Let $\rho$ be a positive measurable (with respect to $d\mu$) function on $U$ that is bounded away from $0$ on every compact subset of $U$. 
It is called a weight.
One can consider the following two $L^2$-norms on sections $\omega\in H^0(U,\eusm K^q)$:
$$
 {\parallel\omega \parallel}' := \sqrt{\int_U \rho\cdot\eta_q(\omega,\bar\omega)d\mu} \quad {\text {and}} \quad{\parallel\omega \parallel}'' :=  \sqrt{\int_U \rho\cdot\Cal H_{U,q}(\omega\otimes \bar\omega)}. 
$$ 
As in \cite{Kol, Chap.\;5.13, Chap.\;5.6-5.8} with obvious modifications, one can compare the  corresponding norms. We have
$$
\qquad \quad\eta_q(\omega,\bar\omega)d\mu = c(n) \eta_q(\xi,\bar \xi) |g|^2 \mu dz_1\wedge \cdots \wedge z_n \wedge d\bar z_1\wedge \cdots\wedge \bar z_n 
$$
and
$$
\quad \Cal H_{U,q}(\omega\otimes \bar\omega)= c(n)\; \varsigma_q |g|^2 dz_1\wedge\cdots \wedge dz_n\wedge d\bar z_1\wedge \cdots \wedge d\bar z_n.
$$
Hence
$$
{{\Cal H_{U,q}(\omega\otimes \bar\omega)} \over {\eta_q(\omega,\bar\omega)d\mu}} = {{\varsigma_q}\over{\eta_q(\xi,\bar \xi)\mu}} \qquad ({\text{on}} \; U)
$$
where $ \epsilon < \eta_q(\xi,\bar \xi)\mu < C $ on $U$ (with constants $\epsilon, C  >0$) because $X$ is compact. We observe that the right-hand side of the last equality is independent of $\omega$.
Thus, if ${\parallel\omega \parallel}' <\infty$ then ${\parallel\omega \parallel}'' <\infty$ because $d\mu$ and $\eta_q$ are $\pi_1(X)$-invariant and $\varsigma_q$ is bounded.



\smallskip

(2.3.3.4) {\it Classical example}\/ I\;(see, e.g., \cite{Kr, Chap.\;III, Sect.\;2-4}). Now, let $X=C$ be a compact Riemann surface (algebraic curve) of $g(C)\geq 2 $.
 Let $U:=\Delta \subset \bold C$ be a disk which is the universal covering of $C$. We set $\lambda_\Delta(z) :=(1-|z|^2)^{-1}$. 
As in (2.3.1.1), we consider a tower of Riemann surfaces with the universal covering $\Delta$:
$$
C=C_0 \leftarrow C_1 \leftarrow C_2 \leftarrow \cdots \leftarrow \Delta; \qquad \bigcap_i Gal(\Delta/C_i) = \{1\} \; (0\leq i<\infty).
$$
In the local coordinate $z$, let $\omega_1 :=g_1(z) (dz)^q$ and $\omega_2 :=g_2(z) (dz)^q$ be weight  $q$ holomorphic forms on $U$ or on the interior of $\Cal D_0 (\Gamma_{\infty i})$. 
We  set  
$$
\Cal H_{\Delta,q}\!: \eusm K_\Delta^q\otimes \bar \eusm K_\Delta^q \rightarrow \Omega_\Delta^{(1,1)}\simeq\eusm K_\Delta\otimes \bar \eusm K_\Delta; \; \omega_1\otimes \bar\omega_2 \mapsto {\sqrt{-1}\over 2}\lambda_\Delta(z)^{2-2q}g_1(z) \overline{g_2(z)}dz\wedge d\bar z. 
$$

In view of the correspondence between automorphic forms on $\Delta$ and the differential forms on the Riemann surfaces, 
we get a compatible sequence of  homomorphisms $\Cal H_{.,q}$ with $\varsigma_q(z) =\lambda_\Delta(z)^{2-2q}$. Recall that $dv_P :=\lambda_\Delta(z)^2 dx\wedge dy$ is the Poincar\'e volume form on $\Delta$  invariant under all holomorphic automorphisms of $\Delta$.
Let $h_{\eusm K^q_\Delta}:= \Cal H_{\Delta,q}/dv_P$ be the corresponding 
invariant Hermitian metric on $\eusm K^q_\Delta $.
\medskip

(2.3.3.5) We will establish a version of a statement
attributed  to Kazhdan by Yau  \cite{Y1, p.\;139}. For a survey of known results
and historical remarks, see a recent article by Ohsawa \cite{O, Sect.\;5}.  In case $U_X$ is the disk $\Delta$, the first proof was given by Rhodes \cite{R}. Recently, McMullen has given a short proof for the disk \cite{M, Appendix}.
\smallskip

 \proclaim{Proposition 1}  With the above notation, we assume that $U$ and all $X_i$'s have the $q$-Bergman metrics for an integer $q$, and the $\Cal H_{.,q}$'s,
defined with a help of the above $\varsigma_q$, are compatible with the tower $(2.3.1.1)$. Then the
$q$-Bergman metric on $U$ equals the limit of pullback on $U$ of the $q$-Bergman metric on $X_i \;(\forall i)$. Furthermore, let  $V$ be the completion of the Euclidean space $E:=\cup V_i$. Then there is a natural isomorphism of projective spaces $\bold P(H^*)\simeq \bold P(V^*)$.
\endproclaim

\demo{Proof} Let $b_{U,q}$ denote the
$q$-Bergman metric on $U$. Set $\tilde b_{U,q}:= \lim \sup  b_i $, where $b_i $ is the pullback on $U$ of the $q$-Bergman metric on $X_i$.  

  First, we establish the inequality $\tilde b_{U,q} \leq b_{U,q}$.
We consider an open precompact exhaustion of $U$, namely: $\{U_\nu \subset U,\; \nu=0, 1, 2,\dots\}$, where $U =\cup U_\nu$,  each $\bar U_\nu$ is compact and $\bar U_\nu\subset U_{\nu+1}$. 
We can take $U_\nu$ to be the interior of the Dirichlet fundamental domain of $\Gamma_{\infty\nu}$ centered at a point $\bold Q$.
 Let $b(\nu)$ denote the $q$-Bergman metric on $U_\nu$. 
 It is well known that  $b_{U,q} = \lim_{\nu\to\infty} b(\nu)$.

Given $U_\nu$, the restriction of $\tau_i$  on $U_\nu$ is one-to-one for $i\gg0$, since $\pi_1(X)$ is residually finite.
 Further,  $b(\nu) > b_i|U_\nu$ for all $i> i(\nu)$. This  establishes that $ b_{U,q} \geq \tilde b_{U,q}$. So, we have  $\tilde b_{U,q}=\lim b_i$. 

Given the metric $\tilde b_{U,q}$,  one can define a natural map
$$
\Upsilon_V: U \longrightarrow \bold P(V^*)
$$
by assigning to each point $u\in U$ the hyperplane in $V$ consisting of elements in $V$ vanishing
at $u$ (compare \cite{Kob2, Chap.\;4.10, p.\;228}). We observe the Cauchy sequence $\{w_k\}\subset E$ vanishes at $u$ if $\lim_{k \to \infty} \omega_k(u) =0$.
The map $\Upsilon_V$ arises from the map of the fundamental domain $U_\nu$ into the corresponding $\bold P(V_\nu^*)$ (for $\nu=0, 1, \dots)$. 
We consider $\bold P(V^*)$ with its Fubini-Study metric. Since $\eusm K_{X_\nu}^q$ is very ample for every $\nu$, $ \Upsilon_V$ is an  embedding.
We get a natural linear map $\delta: H \rightarrow V$.

By assumption, we have an embedding 
$
\Upsilon_H: U \hookrightarrow \bold P(H^*)$. 
Furthermore, we have a natural  embedding $\Upsilon_{VH}: \bold P(V^*) \hookrightarrow \bold P(H^*) $.
 Indeed, a form $\omega \in H$ produces a form on each $X_i$ as follows. First, we restrict $\omega$ on the interior of the fundamental domain of $X_i$. Next, we obtain a measurable form on $X_i$. To obtain a holomorphic form, we apply the Bergman projection. 
We obtain a Cauchy sequence in $E= \cup V_i$ because $\omega$ was square-integrable
on $U$.
We obtain a nontrivial element in $ V$ provided $\omega \not =0$. It remains to show
 that the above construction gives all elements of $V$. 

Let $\bold Q \in U$ be an arbitrary point. An element $v \in V$ is a Cauchy sequence $\{\omega_k\}$ in $E$. We lift each $\omega_k$ on the interior of the corresponding fundamental domain $\Cal D_\bold Q(\Gamma_{\infty\cdot}) \subset U$ and, then, extend by zero outside the interior  of $\Cal D_\bold Q(\Gamma_{\infty\cdot})$. 
We obtain a measurable square-integrable form on $U$. By assumption, $H$ has a reproducing kernel. We apply the Bergman projection \cite{Kol, Chap.\;7.1-7.6} and obtain a holomorphic square-integrable form $\omega^U_k$ in $H$. We get a Cauchy sequence $\{\omega^U_k\}$ in $H$ hence an element in $H$. So, the map $\delta$ is a natural surjection. 

Finally $\Upsilon_H =\Upsilon_{VH}\Upsilon_V.$
 Since the image of $\Upsilon_H$ does not belong to a proper subspace of $\bold P(H^*)$, we obtain the proposition.
\enddemo

\bigskip

In (2.4)-(2.6) below, we assume that $U_X$ is equipped with a $\pi_1(X)$-invariant real analytic Kahler metric $\Lambda_\Cal L$. We consider the induced Riemannian metric on $U_X$ with the volume form $dv_\Lambda$. In Section 3, we shall recall a construction of $\Lambda_\Cal L$ from \cite{T2}.
\smallskip

(2.4) {\it Heat kernel}\/ (see, e.g., \cite{Grig}). 
We consider the heat kernel $p(s,x,y)$ on $U_X$ where 
$
(s,x,y) \in (0,\infty) \times U_X\times U_X.
$
 As a function of $s$ and $x$, the function $p(s,x,y)$ is the smallest positive fundamental solution of the heat equation ${\partial p/\partial s} =\Delta_x p$ where $\Delta_x$ is th Laplace-Beltrami operator on $U_X$ with its  Riemannian metric. 
It is known that on $U_X$ we have (see, e.g., \cite{Grig} and references therein):
\roster
 \item"{(1)}" $p(s,x,y)$ is a $C^\infty$ function of all three variables $(s,x,y) \in (0,\infty) \times U_X\times U_X$
\item"{(2)}" $p(s,x,y) =p(s,y,x)$  $ \quad $ $ \quad $ symmetry
\item"{(3)}" $p(s,x,y) > 0$  $ \qquad $ $\qquad $ $ \quad $ positivity
\item"{(4)}" $\int_{U_X} p(s,x,y) dv_\Lambda (x) = 1$ $ \quad $ stochastic completeness (for an arbitrary Riemannian manifolds $M$, we have $\int_{M} p(s,x,y) dv(x) \leq 1$)
\item"{(5)}" let  $\{U_\nu\}_{\nu \in \bold N} \subset  U_X$ be a precompact open exhaustion 
with smooth boundaries, i.e., $\bar U_\nu\backslash U_\nu$ are smooth and each $\bar U_\nu$ is compact;
let $p_{U_\nu}$ be the smallest positive fundamental solution of the heat equation on $U_\nu$ and set $p_{U_\nu}|U_X \backslash \bar U_\nu:=0$; 
 then  
$$
p_{U_\nu} \leq p_{U_{\nu+1}} \qquad {\text{and}} \qquad p(s,x,y)=\lim_{\nu\to \infty} p_{U_\nu}.
$$
\endroster

When $\bar U_\nu \subset U_X$ is a fundamental domain of a covering of $X$ then $\bar U_\nu \backslash U_\nu$ is not smooth. However, one can approximate the fundamental domains by precompact open regions with smooth boundaries.  
\smallskip

The presence of the Gaussian exponential term in the heat kernel estimates is one of the properties of the heat kernel which little depends on the structure of the manifold in question and reflects the structure of the heat equation.

The heat kernel of the disk $\Delta$ is
decreasing fast as we approach its boundary in $\bold C$. Let $C$ be a compacts Riemann surface of genus at least two with the universal covering $\Delta$.
 Let
 $p_\Delta(s,x,y)$ be the heat kernel on $\Delta$. Set $\bold Q:=0 \in \Delta$. 
We have
$$
\lim_{x \to \bar\Delta\backslash \Delta} p_\Delta(s,x,\bold Q) \lambda^m_{\Delta} (x) =0 
\qquad (\forall m\in \bold Z). \tag{2.4.1}
$$

\smallskip

(2.4.2) {\it Classical example}\/ II (see \!(2.3.3.4)). Set $\bold Q\!:=0 \in \Delta$ and 
$\bold Q_i\!: = \tau_i(\bold Q) \in C_i$.
We consider three volume forms on $\Delta$, namely: the Poincar\'e volume
form $dv_P$,  the Euclidean volume form 
$d\nu:={\sqrt{-1}\over 2}dz\wedge d\bar z =dx\wedge dy$,  and  the following volume form
$$
dv_\Sigma := p^2_\Delta(s,x,\bold Q) dv_P.
$$

It is known that the Poincar\'e metric on $\Delta$ equals the limit of pullback of the Bergman
metric (with respect to the Euclidean volume form $d\nu_i$) on $C_i \; (\forall i)$.
\smallskip

We consider the standard $\pi_1(C_i)$-invariant Hermitian metric $h_{\eusm K_\Delta}$ (see (2.3.3.4)).
Let $H_{p^2}$ be the Hilbert space of square-integrable (with respect
to $dv_\Sigma$ and $h_{\eusm K_\Delta} $) holomorphic sections of $\eusm K_\Delta$.
Also,  $H_{p^2}$ is the space of square-integrable (with respect
to $dv_P$ and the Hermitian metric $\tilde h_{\eusm K_\Delta}:= 
p_\Delta^2(s,x,\bold Q)  h_{\eusm K_\Delta}$) holomorphic sections: 
$$
H_{p^2}:= \biggl\{ \omega \in H^0(\Delta,\eusm K) \;\biggl |\; \Vert \omega  \Vert^2:=  \int_{\Delta} \tilde h_{\eusm K_\Delta}(\omega,\bar \omega) dv_P < \infty   \biggl\}. 
$$
This space is not trivial and has a reproducing kernel.  We obtain a Bergman-type metric on
$\Delta$, denoted by $\beta_{\Delta,p}$. Similarly, for each $i\gg0$, we consider the  Hilbert space $H_i$. Here we keep the same volume form $dv_P$ and consider the
 Hermitian metric $\tilde h_i: = p^2_{C_i}(s,x,\bold Q_i) h_{\eusm K_\Delta}$ where
$p_{C_i}(s,x,y) $ is the heat kernel  of $C_i$.

We obtain
a metric on $C_i \; (\forall i\gg0)$. This metric arises from a natural embedding into a finite-dimensional
projective space $\bold P(H^*_i)$. As in Proposition 1, the limit of pullback on $\Delta$ of the metric
on $C_i \; (\forall i)$ equals the metric $\beta_{\Delta,p}$ because $\pi_1(C)$ is residually finite; 
also see Proposition 1$^\prime$ in Section 5. 

\smallskip

Recall that $\bold B_\Delta(z,\bar z)={1\over \pi}\lambda_\Delta(z,\bar z)^2$, and
$
 \bold P_{\Delta,\bold Q} (\bold p(z,\bar z))=\log\bold B_\Delta(z,\bar z)-\log {1\over \pi}
$ as in (2.2.2),
 where 
 $\bold P_{\Delta,\bold Q}$ is the diastasic potential at $\bold Q$ of the Poincar\'e metric and $\bold p(z,\bar z) \in \Delta$; $\bold B_\Delta$ is the Bergman kernel of $\Delta$ and $\bold B_\Delta (0,0)={1\over \pi}$.
It follows
$$
e^{\bold P_{\Delta,\bold Q} (\bold p(z,\bar z))} = \pi \bold B_\Delta(z,\bar z) =\lambda_
 \Delta(z,\bar z)^2
$$
hence $dv_P=\lambda_ \Delta(z,\bar z)^2dx\wedge dy =e^{\bold P_{\Delta,\bold Q} (\bold p(z,\bar z))} dx\wedge dy$. As expected, we have \lq\lq recovered\rq\rq\/ the Poincar\'e volume form from the diastasic potential. 

The natural map 
from $\Delta$ into  $\bold P(H_{p^2}^*)$, with its Fubini-Study
metric, followed by the Calabi 
{\it flattening out}\/ \cite{C, Chap.\;4, p.\;17} and, then, the projection 
into $\bold C$ produces a bounded domain in $\bold C$ because $zdz \in H_{p^2}$.



\medskip

(2.5) {\it Plurisubharmonic}\/ (psh) and {\it pluriharmonic}\/ (ph) functions on $U_X$. A function $f$ defined in a neighborhood of a point $p \in U_X$ will be psh (respectively ph) if and only if the following holds. For an arbitrary tangent vector $\bold v$ to $U_X$ at $p$, we consider
 $\tau_i( p)$ and $\tau_i(\bold v)$ on $X_i\subset \bold P^{r_i} \; (i\gg 0)$. We assume the latter embedding is nondegenerate and $r_i \gg 0$.
We take a general curvilinear section $C_i\subset X_i$ tangent to $\tau_i(\bold v)$ at $\tau_i( p)$.
 
By our assumption, $C_i$ will be a nonsingular connected curve of genus at least 2.
Indeed, the corresponding linear system has no fixed components because $\pi_1(X)$ is large, in particular, $X$ contains no rational curves. Moreover, the linear system is not composite with a pencil. Hence we can apply Bertini's theorem to obtain a nonsingular connected curve. 
 We, then, consider a connected open Riemann surface $\tau_i^{-1}(C_i) \subset U_X$ (Campana-Deligne theorem \cite{Kol, Theorem 2.14}). Thus, the function $f$ is psh (ph)  at $p$ if and only if its restriction on $\tau_i^{-1}(C_i)$ is subharmonic (harmonic) at $p$.

Let $u$ be a pluriharmonic function on $U_X$. Since $U_X$ is simply connected, there is a holomorphic function $f=u+\sqrt{-1}\tilde u$ on $U_X$ where $\tilde u$ is also pluriharmonic \cite{FG, Chap.\;VI, p.\;318}.
\medskip

(2.6) {\it An easy generalization of some results of Lyons-Sullivan}\/ {\cite{LS, Theorem 3$'$} {\it and Toledo}\/ \cite{To, Lemma 1}. Now, we assume $\pi_1(X)$
is nonamenable.
It follows from \cite{To, Lemma 1} that the space of bounded pluriharmonic functions on $U_X$ is infinite dimensional. In his argument, we replace harmonic functions by pluriharmonic functions.  As in \cite{To, Lemma 1}, the important step is the construction of the map 
$$
\tilde\varphi: L^\infty(U_X) \longrightarrow \text{bounded harmonic functions on}\;  U_X 
$$
by Lyons-Sullivan \cite{LS, Theorem 3$'$, p.\;311}. In view of (2.5), $\tilde\varphi$ produces pluriharmonic functions on $U_X$. Also, the proof of Toledo's lemma \cite{To, Lemma 1} shows that  the linear span of bounded positive pluriharmonic functions is infinite dimensional as well.

\bigskip
\head
3.  Metric $\Lambda_\Cal L$
\endhead

The metric $ \Lambda_\Cal L$ was suggested by a problem of Yau \cite{Y1, Sect.\;6, p.\;139} who proposed to study $\lim_{t \to \infty}{1\over t} g_{X,\eusm K_X^t}$ when $\eusm K_X$ is the ample canonical bundle on $X.$ We do not exclude the case when $\dim X \geq 2$ and $\pi_1(X)$ is Abelian.

The metric $\Lambda_\Cal L$ on $U_X$ is a generalization of the classical Poincar\'e metric though it is not necessary a Bergman-type metric if $\dim X \geq 2$. 
It will depend on the fixed very ample bundle $\Cal L_X$ defining the embedding $\phi: X \hookrightarrow  \bold P^{r}$.
We will define a real analytic potential at every point of $U_X$.
\medskip
(3.1)
First, we will consider the case: $C=X\hookrightarrow \bold P^r$ were $C$ is a connected nonsingular projective curve of genus $g(C)\geq 2$. We will assume the embedding is given by a very ample line bundle $\Cal L_C$ such that 
$$
\Cal L_C \subset \eusm K^{ m}_C,
$$
where $\eusm K_C$ is the canonical bundle and  $ m$ is a suitable  integer.
  We get Bergman-type metrics on $C$ corresponding to $\Cal L_C$ and $\eusm K^{ m}_C $ (see  \cite{Y1, Sect.\;6, p.\;138} and  \cite{Ti, p.\;99}). Also we consider the Poincar\'e metric on $\Delta$. Since $\Delta$ is homogeneous,
$$
\bold B_{\Delta,\eusm K^t}(z,\zeta)= c(t)\bold B^t_{\Delta,\eusm K}(z,\zeta) \qquad  
(\bold B_{\Delta,\eusm K}(z,\zeta)=\pi^{-1}(1-z\bar\zeta)^{-2}),
\tag{3.1.1}
$$ 
where $t\gg 0$ is an integer, $c(t)$ is a known constant  depending on $t$ only, and $\bold B_{\Delta, \eusm K^t}$ denotes the $t$-Bergman kernel  (see, e.g., \cite{FK, p.\;9}, \cite{Kol, (7.7.1)}). It follows
$$
\lim_{t \to \infty}{1\over t} g_{\Delta,\eusm K^t}=\biggl(\lim_{t \to \infty}{1\over t} {\partial^2\log\bold B_{\Delta,\eusm K^t}(z,z) \over \partial z \partial \bar z}\biggl)dz d\bar z =g_{\Delta,\eusm K}.\tag{3.1.2}
$$ 
\medskip

(3.2) Let  $C$ be a sufficiently general nonsingular connected curvilinear section of  $X$\;($g(C)\geq 2$). We consider the inverse image of $C$ on $U_X$. By the Campana-Deligne theorem \cite{Kol, Theorem 2.14}, we obtain a connected open Riemann surface $R=R_C \subset U_X$ in place of the disk $\Delta$. We would like to construct a metric on $U_X$ whose restriction on $R$ is well understood.
Set $\Gamma:= Gal(\Delta /R)$. Let 
$$
\bar \Cal F:= \big\{z\in \Delta \big|\, |Jac_\gamma(z)| \leq 1, \gamma\in \Gamma \big\}, \qquad \big\{\Cal F:= \{z\in \bar\Cal F \big|\, |Jac_\gamma(z)| < 1, \forall\gamma \not = 1\big\}
$$ 
be a fundamental domain of $R$ and the interior of the fundamental domain. 
\smallskip

(3.2.1) High powers of $\Cal L_C:=\Cal L_X |_C$ are squeezed between powers of the canonical bundle on $C$.
 

 For $t\geq1$, let $b_{R,t}$ denote the corresponding Bergman-type metric on $R$
with $\Cal H_{.,t}$ as in the classical example (2.3.3.4). 
Since $R$ is an open Riemann surface, $\Cal L_R$ and $\eusm K_R$ are free bundles.
We do not assume $\eusm K_C$ is very ample.


Let $D_{b,R,t}$ denote the functional element of  diastasis of $b_{R,t}$ at an arbitrary point of $R$.
We lift $b_{R,t}$ and $D_{b,R,t}$ on $\Delta$ and  get
$D_{b,R,t} \leq D_{b,\Delta,t}$ (locally at an arbitrary point of $\Delta$) by Proposition 1. Furthermore, it follows the convergence of the corresponding  holomorphic functions on $R\times \bar R$.
Hence on $\Delta$:
 $$
\lim_{t\to\infty}{1\over t} D_{b,R,t} \leq  D_{b,\Delta,1}=D_{b,\Delta}.
$$
 
Set $b_{R}:=\lim _{t \to \infty}{1\over t}b_{R,t}.$ 
 It will be a  real analytic $Gal(R/C)$-invariant 
metric on $R$. 
We get the metric $b_{R}$ whose diastasic potential $\bold P_{b,R}:=\bold P_{b,R,\bold a}$, where $\bold a \in R$ is the image of the origin $0\in \Delta$, has the prolongation over $R$, i.e., $\bold P_{b,R}$ is a function on $R$. 

One can replace $b_{R,1}$ by $b_{R,m}$, where $m$ is a sufficiently large {\it fixed}\/ number, and repeat the previous argument with $b_{R,m,t}$ in place of $b_{R,t}$. As before, we set 
$$
b_{R,m}:= \lim _{t \to \infty}{1\over t}b_{R,m,t}.
$$ 

We observe that $\Cal L_R^{tm_i} \subseteq \eusm K_R^{mtm_i}$, where $\{m_i\}_{i\in \bold N}$ is a nondecreasing sequence of positive integers. Further, 
$$
b_{R,m} =\lim_{{t,i \to \infty}}{1\over {tm_i}} b_{R,m,tm_i} \qquad {\text {and}} \qquad b_R = \lim_{t,i,m \to \infty} {1\over {mtm_i}} b_{R,m,tm_i}.
$$

We will denote a functional element of the diastasis of $b_{R,m,t}$ by $D_{b,R,m,t}$. As before, the diastasic potentials of $b_{R,m,t}$ and $b_{R,m}$ are functions on $R$.
\smallskip

(3.2.2) Now, let $\Cal L _R$ denote the inverse image of $\Cal L_C$ on $R$. For all $t\gg 0$, we have $\eusm K_R \subset \Cal L^t_R$. Let $H^0_{(2)}(R,\eusm K)$ be the $1$-Bergman space, considered in (3.2.1), giving the embedding
$$
\iota_\eusm K: R \hookrightarrow \bold P ([H^0_{(2)}(R,\eusm K)]^*) \quad (\varsigma_1= 1).
$$

Let $H_t: = H^0_{(2)}(R,\Cal L^t) \;(t\gg 0)$ be the Hilbert space of square-integrable (with respect to the Hermitian metric $h^t$ on $\Cal L^t$ and the restriction of $dv_g$ on $R$) holomorphic global sections. Every $\omega$ in $H^0_{(2)}(R,\eusm K)$ belongs to $H_t$ as in (2.3.3.3) ($\varsigma_1=1$ in the classical example (2.3.3.4)). 
Hence $H_t \not = 0$. We observe that $\omega$ can be viewed as a function on $R$ because $\eusm K_R$ is trivial.

As with the Bergman kernel, we obtain a natural map employing the reproducing kernel of $H_t$:
$$
\iota_{\Cal L^t} : R \rightarrow \bold P (H_t^*).
$$
The map $\iota_\eusm K$ separates points and tangents on $R$. The latter means that, for every point $p\in R$ and the tangent vector $\bold v$ on $R$ at $p$, there exists a section $\omega\in H^0_{(2)}(R,\eusm K)$ such that $\omega(p)=0$ and ${\partial\omega \over \partial\bold v}(p) \not =0$.
 Hence $\iota_{\Cal L^t}$ also separates points and tangents, i.e., $\iota_{\Cal L^t}$ is an embedding as well.

We denote by $g_{R,t}$ the metric on $R$ induced from the Fubini-Study metric on $\bold P(H_t^*)$ via  the map  $\iota_{\Cal L^t}$. We denote by ${1\over t}g_{R,t}$ the inverse image of ${1\over t}$-multiple of the Fubini-Study metric on $\bold P (H_t^*) \; (t\gg 0)$.
Formally, we define
$$
\Lambda_R:=\lim _{t \to \infty}{1\over t}g_{R,t}. 
$$

We will  show the limit exists and $\Lambda_R$ is a real analytic $Gal(R/C)$-invariant Kahler metric on $R$. Let  $D_{R,t}$ denote the functional element of  diastasis of $g_{R,t}$ at an arbitrary point of $R\; (t\gg 0)$. 

Since $\Cal L_R^t \subset \eusm K_R^{mt}$ for all $t\gg 0$ and the {\it fixed
  large}\/ integer $m$, we will be able to show that $D_{R,t}$ is bounded by the corresponding $D_{b,R,m,t}$. 
We will compare the reproducing kernel $B_t$ of $H_t$ and the $tm$-Bergman kernel $\bold B_{mt}$ of $H^0_{(2)}(R, \eusm K^{mt})$.

Let $0\in \bar \Delta_1 \subsetneq \Delta$, where $\bar \Delta_1$ is a closed disk with center $0$ whose radius is close to $1$. Let $\kappa : \Cal O_R \hookrightarrow
\eusm K_R$ be a natural inclusion. Set $\xi:=\kappa(1)$.

On $\Delta \backslash \bar \Delta_1$, we have the estimate
$$
\lambda^{2-2mt} <h^t(\xi, \bar \xi)\!\cdot\!V_g  \tag{3.2.2.1}
$$
provided our $m\gg 0$ (see (2.3.1) and (2.3.3.3)-(2.3.3.4)). 
Indeed,  $\lambda ^{-2m}$ can be made arbitrary small on $ \Delta \backslash \bar \Delta_1 $ by choosing  a sufficiently large $m=m(\Cal L_C)$ (depending on the radius of $\bar \Delta_1 $); recall that $h$ and $V_g$ depend on $\Cal L_X$ only. Furthermore, $h$ and $V_g$ are bounded from below and from above by positive constants because the curve $C$ is compact (see (2.3.3.3)).

By (2.3.3.3), if $\omega \in H_t$ then $\omega \in H^0_{(2)}(R, \eusm K^{mt})$. 
Hence we have a natural inclusion of {\it linear spaces}: 
$$
H_t \subset H^0_{(2)}(R,\eusm K^{mt}).
$$
Since $\eusm K_R^{mt}$ is a free bundle, $\bold B_{mt}$ and $B_t$ can be viewed as functions on $R$.



The Poincar\'e series $\sum_{\gamma \in \Gamma}|Jac_\gamma(z)|^2$ is uniformly convergent. 
As well known, this yields that
 every compact subset of $\Delta$, in particular $\bar \Delta_1$, is covered by a finite number of $\gamma \bar\Cal F \; (\gamma \in \Gamma)$, where $\bar\Cal F$ is the fundamental domains defined in (3.2).

Similarly, for each $i$, set $C_i:=\tau_{i0}^*(C)$ where $C\subset X\subset \bold P^r$ is a general curvilinear section. Let $Gal(R/C_i)$ be the corresponding Galois group, and let $\bar\Cal F_i \subset R$ be the 
fundamental domain of $Gal(R/C_i)$. Relative Poincar\'e series are discussed, e.g, in \cite{Dr}. As above, for each $i$, 
every compact subset of $R$ is covered by a finite number of $\gamma \bar\Cal F_i \subset R\; (\gamma \in Gal(R/C_i))$. Also, the estimate (3.2.2.1) holds on $R\backslash\Psi$ where $\Psi \subset R$ is a compact subset.

Now, we consider an arbitrary $\omega  \in H_t$. Let $\parallel\cdot\parallel'$ and $\parallel\cdot\parallel''$ denote the norms in $H_t$ and $H^0_{(2)}(R,\eusm K^{mt})$, respectively.
Then $ \parallel  \omega \parallel'' \leq \parallel  \omega \parallel'  
$ for all $\omega \in H_t$ by (2.3.3.3) and the above estimate (3.2.2.1). 
 Indeed, the inequality (3.2.2.1) can be violated on at most a finite number of fundamental domains of $C$, and $\Delta$ (and $R$) are covered by {\it infinitely}\/ many corresponding fundamental domains.
 
It follows  $B_t \leq \bold B_{mt}$ on $R$ (see, e.g., \cite{FK, p.\;6}, \cite{Sha, 
Chap. 6.17, Sect.\;51}), as well as $B_t  \ll  \bold B_{mt}$ on $R$ by Aronszain
(see (2.3.2.3)). 
As in  (2.2.2), 
 $$
\lim_{t\to \infty}{1\over t}D_{R,t}    \leq    \lim_{t\to \infty}{1\over t}D_{b,R,m,t}. 
\tag{3.2.2.2}  
$$

It follows the uniform convergence of the corresponding  
holomorphic functions on $R\times \bar R$.  Furthermore, ${1\over t}D_{R,t}$ generates a global function on $R$. It follows the diastasic potential of the metric $\Lambda_R$ on $R$, $\bold P_R:=\bold P_{R,\bold a}$  where $\bold a \in R $ is the image of the origin $0\in \Delta$, has the prolongation over $R$. 

Finally, $g_{R,t}$ is $Gal(R/C)$-invariant as, e.g., in \cite{FK, pp.\;10-11 or pp.\;188-191},
\cite{Kob1,Theorem 3.1}, \cite{Sha, Chap.\;6.17, Sect.\;52, Theorem 1}.

\proclaim{Proposition-Definition 2} With the notation of $(3.2)$, let
$\Lambda_R :=\lim_{t \to \infty}{1\over t} g_{R,t}.$
Then $\Lambda_R$
 is a real analytic $Gal(R/C)$-invariant Kahler metric. The diastasic potential\/
 $\bold P_R:=\bold P_{R,\bold a}$, where $\bold a\in R$ is the image of the origin\/ $0\in \Delta$, has the prolongation over $R$. 
\endproclaim

\medskip

(3.3) Now, we return to the situation in (2.3.1) with $U=U_X$. For each positive integer $t$, the Hermitian metric $h$ on $\Cal L_X$ induces a Hermitian metric $h^{t}$ on $\Cal L_X^{t}$ as well as on all inverse images of  
$\Cal L_X^{t}$ on the coverings of $X$. 

\smallskip
(3.3.1)
	We choose an orthonormal basis $(s^{t}_0, \dots, s^{t}_{r_{t}})$ of   $H^0(X, \Cal L^{t}_X)$ with respect to $dv_g$ and $h^{t}$.  
We have an inner product and a natural embedding:
$$
\langle s^{t}_\alpha, s^{t}_\beta \rangle := 
\int_X h^{t}(s^{t}_\alpha, s^{t}_\beta)dv_g; \quad  \phi_{X,{t}}: X\hookrightarrow \bold  P^{r_{t}}:=\bold  P(H^0(X, \Cal L^{t}_X)^*).
$$
Let $g_{FS}$ denote the corresponding standard Fubini-Study metric on the projective space. As in Yau \cite{Y1, Sect.\; 6, p.\; 139} (see also Tian  \cite{Ti}), the ${1\over {t}}$-multiple of  $g_{FS}$ on $\bold P^{r_{t}}$ restricts to a Kahler metric  on  $X$:
$$
g_{X,{t}}:={1\over {t}} \phi^*_{X,{t}}g_{FS}.
$$
\smallskip

(3.3.2) We consider the finite coverings of $X$. The bundles $\tau_{i0}^*\Cal L_X$  $(1\leq i<\infty)$ are ample by the Nakai-Moishezon ampleness criterion. However, $\tau_{i0}^*\Cal L_X$'s are not necessary very ample  bundles.  
Let $\{m_i\} \;(i\geq 1)$ be a nondecreasing sequence of positive integers such that the bundle $(\tau_{i0}^*\Cal L_X)^{m_i}$ is very ample. Then the bundle $(\tau_{i0}^*\Cal L_X)^{tm_i}$  defines a natural embedding
$$
\phi_{X_{i,tm_i}}: X_i \hookrightarrow \bold P^{r_{tm_i}}
$$ 
into an appropriate projective space.
We get a metric  $g_{X_{i,tm_i}}:= {1\over {tm_i}} \phi^*_{X_{i,tm_i}}g_{FS}$  on $X_i$ and the corresponding diastasic potential.  

For all $i=1,2, \dots$, we consider the metrics $g_{X_{i,tm_i}}$  as $t \rightarrow \infty$.
\smallskip

(3.3.3)  We consider pullbacks on $U_X$  of the  metrics $g_{X_{i,tm_i}}$  and the corresponding diastasises. We will establish that the functional elements of the diastasises converge at a point $\bold p \in U_X$, and we will obtain a real analytic strictly plurisubharmonic functional element at $\bold p$. For points of $U_X$, such functional elements will define the desired Kahler metric $\Lambda_\Cal L$ on $U_X$. 

Let $H_\infty$ denote the hyperplane at infinity in $\bold P^r$.  We can and will assume that  the functional elements of diastasises generate functions on the preimages of $X\backslash H_\infty$ on $X_i$'s (see \cite{C, Chap.\;4}).

	We assume the point $\bold p$ does not lie at infinity.  We consider a small compact neighborhood $G\subset U_X$ of $\bold p$. The pullbacks  on $U_X$ of the above functions produce functions on $G$. 

First, we will establish the pointwise convergence of the pullbacks of this functions in a small neighborhood of   $\bold p$. 
We will make use of the fundamental property of the diastasis.

We take a general curvilinear section $C\subset X\subset \bold P^r$ whose inverse image on $U_X$ contains $\bold p$. This inverse image will be a connected open Riemann surface by the Campana-Deligne theorem \cite{Kol, Theorem 2.14}.  We, then, apply Proposition-Definition 2. 

The pointwise convergence is independent of the curve $C$ because  we have the metric $g_{X_{i,tm_i}}$ on each $X_i$ $(0\leq i <\infty)$. By the fundamental property of the  diastasis, we get the same functional element of diastasis at $\tau_i(\bold p)\in \tau^{-1}_{i0}(C)$, independent of the choice of the Riemann surface $\tau^{-1}_{i0}(C)$, that is, of the curve $C$.

We obtain a function on  $G\times \bar G$. By Hartogs' theorem (separate analyticity implies  joint analyticity), this function will be holomorphic. 
In particular, it follows  the functional elements of diastasises converge uniformly on $G$ by  Dini's monotone convergence theorem, and we get the uniform convergence on $G\times \bar G$ of the corresponding holomorphic functional elements  (the complexifications) to a holomorphic functional element. 

We obtain a {\it real}\/ analytic functional element at the fixed point $\bold p$, {\it denoted}\/ by $D_{U_X}(\bold p,u)$.
\smallskip

 (3.3.4) It is easy to see that $D_{U_X}(\bold p, z(u), \bar z(u))$ is strictly plurisubharmonic, where $z=(z_1,\dots,z_n)$ are coordinates in a neighborhood with origin at $\bold p$. Indeed, we take an arbitrary tangent vector $\bold v$ to $U_X$ at $\bold p$. We consider $\tau_i(\bold p)$ and $\tau_i(\bold v)$ on $X_i\subset \bold P^{r_{m_i}}\;  (i\gg 0)$. We take a general curvilinear section $C_i\subset X_i$ tangent to $\tau_i(\bold v)$ at $\tau_i(\bold p)$. 

By our assumption, $C_i$ will be a nonsingular connected curve of genus at least 2.
Indeed, the corresponding linear system has no fixed components because $\pi_1(X)$ is large, in particular, $X$ contains no rational curves. Moreover, the linear system is not composite with a pencil. Hence we can apply Bertini's theorem to obtain a nonsingular connected curve. 
 We, then, consider a connected open Riemann surface $\tau_i^{-1}(C_i) \subset U_X$. 
We get $D_{U_X}(\bold p,z(u), \bar z(u))$ is strictly plurisubharmonic  at $\bold p$.

Finally, the $n\times n$ matrix $(h_{\alpha\beta})$, where
$$
h_{\alpha\beta}(z(u),\bar z(u)): = {\partial^2 \over \partial z_\alpha \partial \bar z_\beta} D_{U_X}(\bold p,z(u), \bar z(u)),
$$
defines the desired real analytic Kahler metric $\Lambda_\Cal L$. 
By (2.2.1) and the  fundamental property of the diastasis,  $D_{U_X}(\bold p,u)$ is, 
in fact, the {\it diastasic} potential of $\Lambda_\Cal L$ at $\bold p$.

Thus, we have established the following

\proclaim{Proposition-Definition 3} We assume $\pi_1(X)$ is large and residually finite,  and a general curvilinear section $C\subset X$ has $g(C)\geq 2$. Then $U_X$ is equipped with a real analytic $\pi_1(X)$-invariant Kahler metric, denoted by $\Lambda_\Cal L$.
 The restriction of $\Lambda_\Cal L$ on $R_C$, the inverse image on $U_X$ of a general curvilinear section $C\subset X$, is the metric $\Lambda_R$ on $R_C$. 
\endproclaim
\bigskip

\head
4.  Metric $\Sigma_\Cal L$
\endhead 

(4.1) We assume, in addition,   $\pi_1(X)$ is nonamenable. We will construct a metric $\Sigma_\Cal L$ and prove $\eusm K_X$  is ample. Let  
$dv_\Lambda=V_\Lambda \prod^n_{\alpha=1}\bigl({\sqrt{-1}\over2} \cdot\!dz_\alpha\!\wedge\! d\bar {z}_\alpha \bigl)
$ denote the volume form of $\Lambda_\Cal L$,
where $V_\Lambda$ is a locally defined positive function. As before, $p_{U_X}(s,x,y)$ denotes the heat kernel on $U_X$. Let $\bold Q \in U_X$ be a fixed point. 

We consider the following    volume form on $U_X$:
$$
dv_\Sigma : = p^2_{U_X}(s,x,\bold Q)\!\cdot\!dv_\Lambda .
$$

A harmonic function $u$ on the disk $\Delta$ is the real part of one and only one holomorphic function $f_u=u+\sqrt{-1}\tilde u \in \text {Hol}(\Delta)$ such that  $f_u(0)=u(0)$. By a theorem of M. Riesz (communicated to the author by Demailly; see \cite{Ru, 17.24-17.26}) if $u$ is  $L^2$ with respect to the Lebesque measure then
$$
\Vert f_u\Vert \leq A \Vert u\Vert , \tag{4.1.1}
$$
where $\Vert \cdot \Vert$ is the norm in $L^2$ and $A$ is a constant.
\smallskip

(4.2) Now, we consider two general subspaces $\bold P^{n-1}\subset \bold P^r$ and $\bold P^{r-n} \subset \bold P^r$, $ \bold P^{n-1} \cap\bold P^{r-n}=\emptyset$. We obtain a fibration $\Cal F_X$ of curves over an open subset of the $\bold P^{n-1}$ whose general member is a nonsingular curve in $X$ of genus at least 2.
The inverse image of $\Cal F_X$ on $U_X$ produces a family $\Cal F_U$ whose general member is a connected open Riemann surface on $U_X$.

A bounded pluriharmonic function $u$ on $U_X$ will be square-integrable with respect to  $dv_\Sigma$. 
Furthermore, the restriction of $u$ on the general member $R\in \Cal F_U$ will be $L^2$ with respect to the Lebesque measure on $R$.   
By (4.1.1) and the Fubini-Tonelli theorem, the corresponding  holomorphic function
$f=u+\sqrt{-1}\tilde u$ on $U_X$ will be square-integrable with respect to  $dv_\Sigma$.
\smallskip

(4.3) Let $H_\Sigma$ be the Hilbert space of square-integrable on $U_X$ (with respect to $dv_\Sigma$) holomorphic functions on $U_X$:
$$
H_\Sigma := \biggl\{ \varphi \in \text {Hol}(U_X) \;\biggl |\; \Vert \varphi  \Vert^2:=  \int_{U_X}|\varphi (x)|^2 dv_\Sigma < \infty   \biggl\}. 
$$
This Hilbert space is infinite dimensional and has a reproducing kernel, denoted by
$\bold B_\Sigma (z,\bar z)$ (see, e.g., \cite{FK, pp.\;8-10 or pp.\;187-189});
$\bold B_\Sigma (z,\bar z)$ is a function on $U_X$. 
 Set 
$$
 g_{\alpha \beta}:={\partial^2\log\bold B_\Sigma (z,\bar z) \over \partial z_\alpha \partial \bar z_\beta}.  
$$

\smallskip
(4.4) The differential form
$
ds_\Sigma^2 := \sum^n _{\alpha, \beta =1}g_{\alpha \beta} d z_\alpha d \bar z_\beta
$
is called a Bergman form.  Clearly  $ds_\Sigma^2$ is Hermitian.  We claim it is positive definite (meaning $\log \bold B_\Sigma (z,\bar z)$ is strictly plurisubharmonic), i.e., for any vector  $\bold w\in \bold T_{\bold p,U_X}$, $\bold w \not=0$, at an arbitrary point $\bold p\in U_X$:
$$
 \sum^n _{\alpha, \beta =1}g_{\alpha \beta} w_\alpha  \bar w_\beta > 0 \qquad  (\bold w=(w_1, \dots, w_n)).
$$
Now, we fix $\bold p$ and $\bold w\in \bold T_{\bold p,U_X}$. To prove the positivity, we consider the set
$$
\Cal E :=\{ \varphi \in H_\Sigma \; \bigl | \;\varphi (\bold p)=0, \langle \nabla \varphi, \bold w \rangle =1  \bigl\}.
$$
 A priori, it is not obvious that  $\Cal E \not = \emptyset$. First, we assume  $\Cal E \not = \emptyset$ and show that
$$
\min_{\varphi \in \Cal E}  \| \varphi \|^2 = {1\over \bold B_\Sigma (z,\bar z) \sum _{\alpha, \beta }g_{\alpha \beta} w_\alpha  \bar w_\beta}, 
$$
where $\bold B_\Sigma (z,\bar z)$ and $ g_{\alpha \beta} $ are computed at $\bold p$, hence $ \sum _{\alpha, \beta }g_{\alpha \beta} w_\alpha  \bar w_\beta  > 0$.  Our argument is similar to the one in \cite {Sha, Chap.\;6.17, Sect.\; 52}. We will briefly recall his argument.

Let  $\{ \varphi_\sigma\} \subset H_\Sigma$ denote the complete orthonormal system in $H_\Sigma$.  Let  $\varphi = \sum_\sigma a_\sigma \varphi_\sigma$. Then our problem is to find $\min \sum_\sigma|a_\sigma|^2$ under the conditions  at $\bold p$:
$$
\sum_\sigma a_\sigma \varphi_\sigma = 0, \qquad \sum_\sigma a_\sigma \langle \nabla\varphi_\sigma,\bold  w \rangle =1.
$$
We employ the method of Lagrange multipliers. The uniqueness is easy, provided we have a solution (see \cite {Sha, Chap.\;6.17, Sect.\; 51}).

For the extremal values of  $a_\sigma$, we obtain as in \cite {Sha, Chap.\;6.17, Sect.\; 52}:
$$
\sum_\sigma a_\sigma \bar a_\sigma = {1 \over \bold B_\Sigma (z,\bar z)\sum_{\alpha, \beta} g_{\alpha\beta} w_\alpha \bar w_\beta}.
$$
By the above, we get a Bergman-type metric $ds_\Sigma^2$ and a holomorphic immersion into an infinite-dimensional projective space as in \cite{Kob2, Chap.\;4.10} provided  $\Cal E \not = \emptyset$.
\medskip

(4.5)
It remains to verify 
 $\Cal E \not = \emptyset$.   
We consider a {\it general}\/ curvilinear section $C_i$ of $X_i\subset \bold P^{r_i}\; (i\gg 0)$ through the image of $\bold p$ and $\bold Q$ on $X_i$ that is tangent to $\tau_i(\bold w)$ (see (3.3.4)). Let $R:=R_{C_i} \subset U_X$ be the corresponding open connected Riemann surface\; 
($\Delta\rightarrow R \rightarrow C_i$). 

Now, let $u$ be a non-constant bounded positive pluriharmonic function on $U_X$, and $f$ the corresponding holomorphic function on $U_X$ (see (2.5)).

Let $z_1$ be the coordinate in $\Delta\subset \bold C$ with origin at $0\in \Delta$ ($0 \mapsto  \bold p\in R$). Let $f_R$ denote the restriction of $f$ on $R$. We can assume
 $f_R$ is not a constant. By abuse of notation, we denote its pullback on
$\Delta$ by $f_R$ as well. Let
$$
f_R(z_1)= f_R(0) + a_1z_1 + a_2z_1^2 +\cdots
$$ 
be the Taylor expansion around $0\in \Delta$. 
Let $k$ be the first index such that $a_k \not =0$.
We differentiate $f$ in the direction $\bold w$ at $\bold p$ exactly $k$ times. We obtain a holomorphic function $f^{(k)}$ and the corresponding function $f_R^{(k)}$.

We consider the following generalization of the Schwarz-Pick inequality in Dai-Pan \cite{DP, Theorem 1.2}:
$$
|f_R^{(k)}| \leq {{2k!\; \Re{f_R(z_1)}} \over {(1-|z_1|^2)^k}}\;(1+|z_1|)^{k-1}.
$$
It follows the function $f^{(k)}$ belongs to $\text{H}_\Sigma$  by (2.4.1).
\bigskip

(4.6) Thus, the metric $ds^2_\Sigma$ arises via the immersion 
$$
\Upsilon : U_X \longrightarrow \bold P(H_\Sigma^*); \quad \Upsilon^* ds^2_{\bold P(H_\Sigma^*)}=ds^2_\Sigma
$$
 \cite{Kob2, Chap.\;4.10, p.\;228}.
If $\{\varphi_j\}$ is an orthonormal basis of $ H_\Sigma$ then  $\Upsilon$ is
given by $u\mapsto [\varphi_0(u) : \varphi_1(u) :\dots]. $
The fundamental group $\pi_1(X)$ acts on $H_\Sigma$ as follows
$$
T_\gamma :\varphi \mapsto (\varphi\circ\gamma)\cdot Jac_\gamma \qquad (\gamma
\in \pi_1(X))
$$
where $Jac_\gamma$ is the (complex) Jacobian determinant of $T_\gamma$
(see, e.g., \cite{FK, pp.\;10-11 or pp.\;188-191}, \cite{Kob1,Theorem 3.1}, \cite{Sha, 
Chap.\;6.17, Sect.\;52, Theorem 1}). The map $T_\gamma$ is a unitary isomorphism of
$H_\Sigma$. 

We get a $\pi_1(X)$-invariant volume form on $U_X$ defined in local Bochner
canonical coordinates (\cite{Bo, p.\;181}, \cite{C, Chap.\;3, Prop.\;7}) as follows:
$$
\bold B_\Sigma (z,\bar z)\!\prod^n_{\alpha=1}\Bigl({\sqrt{-1}\over2} dz_\alpha\!\wedge  d\bar {z}_\alpha \Bigl) =\bold B_\Sigma (z,\bar z)\!\prod^n_{\alpha=1}\bigl(dx_\alpha\!\wedge dy_\alpha \bigl)
 \; \Bigl(\bold B_\Sigma(z,\bar z)=\sum^\infty_{j=0}|\varphi_j(z)|^2 \Bigl).
$$
This volume form is $\pi_1(X)$-invariant as, e.g., in  \cite{FK, pp.\;10-11 or pp.\;188-191},
\cite{Kob1, Theorem 3.1}, \cite{Sha, Chap.\;6.17, Sect.\;52, Theorem 1}. 

Thus $X$ admits a volume form (and the corresponding Hermitian metric on $\eusm K^{-1}_{U_X}$) with negative Ricci form  \cite{Kob2, Chap.\;2.4.4}. It follows $\eusm K_X$ is ample by Kodaira. 
 
\bigskip

\head
5.  Metric $\beta$ and uniformization
\endhead 

We assume $\pi_1(X)$ is nonamenable. In this section, we fix an integer $q\gg 0$ and assume $\eusm K^q $ is  very ample on every finite covering of $X$  (see (2.3.3)). Now, we  consider the  $\pi_1(X)$-invariant metric $\Lambda_{\eusm K^q}$ on $U_X$. Let $\bold Q \in U_X$ be our fixed point
(see (4.1)). We will consider a Bergman-type  metric with a weight, denoted by $\beta:=\beta_{\eusm K^q}$, similar to the metric $\Sigma_\Cal L$.

By (4.6), we get a $\pi_1(X)$-invariant Hermitian metrics on $\eusm K^{-1}_{U_X}$ and $\eusm K_{U_X}$. The metric on $\eusm K_{U_X}$ is denoted by $\eta_1$.
\smallskip

(5.1) Now, let $\eta_q:= \eta_1^q$ be the  $\pi_1(X)$-invariant Hermitian metric on $\eusm K^q_{U_X}$.
Let $H_\beta$ be the Hilbert space of square-integrable, with respect to the Hermitian 
metric $\tilde\eta_q:= p^2_{U_X}(s, x, \bold Q) \eta_q $  and  $dv_\Lambda$
($\Cal L=\eusm K^q$), 
holomorphic sections of $\eusm K^q_{U_X}$:
$$
H_\beta := \biggl\{ \omega \in H^0(U_X,\eusm K^q) \;\biggl |\; \Vert \omega  \Vert^2:=  \int_{U_X} \tilde\eta_q(\omega,\bar \omega) dv_\Lambda < \infty   \biggl\}. \tag{5.1.1} 
$$
By Section 4, this Hilbert space is infinite dimensional  because of the natural inclusion $O_{U_X} \hookrightarrow \eusm K^q_{U_X}$. It has a reproducing kernel. As in (4.4)-(4.5), we obtain a positive definite
Bergman form $ds^2_\beta$ and a metric  on $U_X$, denoted by $\beta$.

This metric arises from a natural immersion of $U_X$ into the projective space $\bold P(H_\beta^*)$ with its Fubini-Study metric \cite{Kob2, Chap.\;4.10}. We will show this immersion is actually an embedding. This will follow at once from Proposition 1$'$ below whose proof is a trivial generalization of Proposition 1 because $\pi_1(X)$ is residually finite.
\smallskip

(5.2) Let $p_i(s,x,y)$ be the heat kernel on $X_i$.  Set $\bold Q_i:=\tau_i(\bold Q)\in X_i $.  On each $X_i$, we consider the  Hermitian metric
$$
\tilde \eta_{i,q}: = p^2_i(s,x,\bold Q_i) \eta_{q}.
$$
For every $i$, we consider the Hilbert space $H_{i,\beta}$ defined similarly to (5.1.1)
with $\tilde\eta_{i,q}$ in place of $\tilde\eta_q$.
We obtain  a metric  on $X_i$, denoted by $\beta_i :=\beta_{i,\eusm K^q}$.
This metric arises from a natural embedding into a finite-dimensional projective space $\bold P(H_{i,\beta}^*)$. 

\proclaim{Proposition 1$'$} With the above notation, the metric $\beta$ equals the limit of pullback on $U_X$ of the metrics $\beta_i\; (\forall i)$.
\endproclaim
\smallskip

(5.3) {\it End of proof of the theorem.}\/ We assume $U_X$ is equipped with the metric $\beta$. We will apply the Calabi {\it flattening out}\/ $\eurb S_\bold Q$ (a generalized stereographic projection) of $U_X$ from the point 
$\bold Q$
into the unitary space, i.e., the Fubini space  $\bold F(\infty,0)$ \cite{C, Chap.\;4, p.\;17}.
By \cite{C, Chap.\;4, Cor.\;1, p.\;20}, the whole projective space $\bold F(\infty,1)$,
except the antipolar hyperplane $A$ of $\bold Q$, can be flatten out into $\bold F(\infty,0)$.

If $n=1$ then the uniformization theorem is due to Koebe and 
Poincar\'e, independently (we will employ a little more precise statement;  see (2.4.2)). 
In this subsection, if $n=\dim X \geq 2$ then we assume the uniformization 
theorem is valid for all dimensions less than $n$. Precisely, for $k\leq n-1$, the map 
from $U_X\;(\dim U_X=k)$ into the corresponding $\bold P(H^*_\beta)$ followed 
by the Calabi {\it flattening out}\/ and, then, the projection into $\bold C^k$ 
produces a bounded domain in $\bold C^k$. By the Poincar\'e residue map \cite{GH,
Chap.\;1.1, p.\;147}, we can assume that the restriction of $\eusm K^q_{X_i}\; (q\gg
0)$ on the general hyperplane section of $X_i\subset \bold P^{r_{m_i}}$ is a very ample pluricanonical bundle.

Let $(x^\sigma)$ be  arbitrary Bochner canonical coordinates 
 with origin at $\bold Q$  
in the ambient projective space minus the antipolar hyperplane. The canonical coordinates are unique up to homogeneous
unitary transformations \cite{C, Chap.\;3, Prop.\;7, p.\;14}.
We put a new Kahler metric in a {\it neighborhood} of $\bold Q$ in $ \bold F(\infty,0)$
by replacing $\bold P_\bold Q(\bold p)$ with the diastasic potential
$$
\bold P^{f\!l}_\bold Q(\bold p):= e^{\bold P_\bold Q(\bold p)}-1 = \sum_{\sigma=1}^\infty
|x^\sigma(\bold p)|^2.
$$

Recall that $U_X$ is compete. As in \cite  {C, Chap.\;3, Theorem 7, p.\;15}, we can assume
$$
\bold P^{f\!l}_{U_X\!,{\bold Q}}(\bold p) = \sum_{\sigma=1}^n
|x^\sigma(\bold p)|^2 \qquad (\bold p \in U_X). \tag{5.3.1}
$$ 
A priori, the flattening out can be performed only in a neighborhood of $\bold Q$ in $U_X$. However $\bold P_{U_X\!,\bold Q}$ can be prolonged over  $U_X$ by Lemma A in the Appendix.
Thus, there is a Bochner canonical coordinate system in $U_X \subset \bold F(\infty,1)$
with origin at $\bold Q$ that covers the whole $U_X$.
 Hence we can flatten out the whole $U_X$
because the Bochner canonical coordinates can be extended as  meromorphic functions 
\cite {C, Chap.\;4,  Cor.\;3, p.\;20} and, in our case, they will be holomorphic functions.
 The sum in (5.3.1) is from $\sigma=1$ to $\sigma=n$ because of the natural projection 
from $U_X\subset  \bold F(\infty,0)$ into $\bold F(n,0)$. In fact, the 
Bochner canonical coordinates produce a holomorphic bijection of $U_X$ into 
$\bold C^n= \bold F(n,0)$.

It remains to show $\eurb S_\bold Q(U_X)\subset \bold C^n$ is a {\it bounded}\/ domain,
 i.e., $|x^\sigma (\bold p)| <\eurb L$ for a constant $\eurb L$, $1\leq \sigma \leq n$ and $\forall \bold p\in U_X$.
Let $ \Cal C$ be a general hyperplane section of $X\subset \bold P^r$ through $\tau(\bold Q)$.
Set $V_\Cal C^{n-1}:=\tau^*(\Cal C)$. By induction hypothesis, $V_\Cal C^{n-1}$ is a bounded domain for $n\geq3$. It is an open Riemann surface for $n=2$. 

The diastasic potentials $\bold 
P_{V_\Cal C^{n-1},\bold Q}$ and $\bold P^{f\!l}_{V_\Cal C^{n-1},\bold Q}$ are functions. By induction hypothesis,
$\bold P^{f\!l}_{V_\Cal C^{n-1},\bold Q}$ is bounded on $\eurb S(V_\Cal C^{n-1})$ by 
 $\eurb r_{V_\Cal C^{n-1}}$,  where 
$$
\eurb r_{V_\Cal C^{n-1}}:= \sup_{\bold p\in\eurb S(V_\Cal C^{n-1})}\Bigl(\sum_{\sigma=1}^{n-1}|x^\sigma(\bold p)|^2\Bigl)   
$$ 
depends only on  the Kahler metric because the diastasis
is uniquely determined by the Kahler metric. Furthermore, we have assumed that $x^1,\dots, x^{n-1}$ are Bochner canonical coordinates in $V_\Cal C^{n-1}$.
Clearly $r_{V_\Cal C^{n-1}} <\infty$ if and
only if $|x^\sigma(\bold p)|^2$ is bounded for all $\sigma \; (1\leq\sigma \leq n-1)$ and all $ \bold p\in\eurb S(V_\Cal C^{n-1})$.

For $n\geq 3$, $V_\Cal C^{n-1}$ is a bounded domain 
by induction hypothesis. 
For $n=2$, we get $\eurb r_{V_\Cal C^1} \not = \infty$ by the construction of the metric $\Lambda_R$ 
(see the proof of Proposition-Definition 2 in (3.2.2), especially (3.2.2.3)).
The number $\eurb r_{V_\Cal C^{n-1}}$ is independent of the choice of Bochner canonical coordinates by \cite {C, Chap.\;3, Theorem 7, p.\;15}.



Finally, one can choose $\Cal C'\subset X$ (similar to $\Cal C$) of dimension $n-1$ 
such that
$$
\eurb r_{V_{\Cal C'}^{n-1}}:= \sup_{\bold p\in\eurb S(V_{\Cal C'}^{n-1})}\Bigl(\sum_{\sigma=2}^{n}|y^\sigma(\bold p)|^2 \Bigl) < \infty,
$$
where $(y^\sigma)$ are Bochner canonical coordinates with origin at $\bold Q$
in the ambient projective space minus the antipolar hyperplane as well as in $U_X$ 
when $1\!\leq\!\sigma\!\leq\!n$.  Furthermore, we can assume $y^2,\dots,y^n$ are Bochner canonical coordinates in $V_{\Cal C'}^{n-1}$, and we have chosen $\Cal C'
$ such that $y^n=x^n$ on $U_X$. The latter can be achieved (see (2.5) or (3.3.4)).
It follows $|x^n(\bold p)|^2 <\eurb r_{V_{\Cal C'}^{n-1}}$.
Hence $ \eurb S(U_X) $ is a bounded domain.

This proves the theorem.



\bigskip
\head
Appendix
\endhead

We will prove the following theorem conjectured by Shafarevich (1972).

\proclaim{Theorem A} Let $X \hookrightarrow  \bold P^{r}$ be a nonsingular connected projective variety of dimension $n\geq 1$ with {\it large} and {\it residually finite} fundamental group $\pi_1(X)$. Then $U_X$ is a Stein manifold.
\endproclaim

\demo{Proof} We will prove the theorem under an additional assumption that a sufficiently general nonsingular connected curvilinear section $C\subset X$ has the genus $g(C)\geq 2$. Otherwise, $\pi_1(X)$ is  Abelian by the Campana-Deligne theorem \cite{Kol, Theorem\;2.14}, and the conjecture is well known when $\pi_1(X)$ is Abelian or, even, nilpotent (Katzarkov \cite{Ka}). 
\medskip

(A.0)  The idea of proof of the theorem is similar to the one by 
Siegel \cite{S1}. 
He established  that
if  $U$ is a connected bounded domain in $\bold C^n$ covering a compact complex manifold $Y$ then $U$ is a domain of holomorphy. We will sketch his argument.

He considers the Bergman metric on $U$  (see, e.g., \cite{Kob1}, \cite{Kob2, Chap.\;4}).  It is complete since $Y$ is compact. 
Recall the fundamental property of the Bergman metric, namely, it determines a natural  isometric embedding of $U$ into an infinite-dimensional projective space with a Fubini-Study metric. Infinite-dimensional projective spaces were considered by Bochner \cite{Bo, p.\;193}, Calabi  \cite{C, Chap.\;4} and Kobayashi \cite{Kob1, Sect.\;7}. According to Kobayashi \cite{Kob1, p.\;268},  the idea of using square-integrable forms on arbitrary  manifolds can be found in Washnitzer \cite{Wa}.

 Let $\bold B(z,\bar z)$ denote the Bergman kernel of $U$.  Siegel proves that  $\log \bold B(z,\bar z)$
goes to infinity on any infinite discrete subset $T\subset U$.
Hence  $U$ is a domain of holomorphy (equivalently, holomorphically compete or Stein domain)  by Oka's solution of the Levi problem. 

We observe that $Y$ is  a projective variety by the Poincar\'e ampleness theorem \cite{Kol, Theorem 5.22}. In 1950s, Bremermann proved  that an arbitrary bounded domain in $\bold C^n$ with complete Bergman metric is Stein \cite{Kob2, Theorem 4.10.21}. 

Now, Oka's solution of the Levi problem for domains in $\bold C^n$ admits a generalization due to Grauert (manifolds) and Narasimhan (complex spaces) (see \cite{N}). Thus, our aim is to define a metric on the manifold $U_X$ and a strictly plurisubharmonic function on $U_X$ that goes to infinity on an arbitrary infinite discrete subset of $U_X$.  The function generated by the diastasic potentials of the metric will be such a function.

We assume  $U_X$ and all $X_i$'s are equipped with the metric $\Lambda_\Cal L$. Let $I, R\subset U_X$ be two subsets with $I$ compact. Let  $d(u,p)$ denote the distance function on  $U_X$ with its Riemannian structure induced by $\Lambda_\Cal L$.
 We set
$$
  \xi(I,R):=\sup_{u\in I} [\inf_{p\in R} d(u,p)].
$$

 We say that a sequence of subsets $\{R_\gamma\}_{\gamma\in \bold N}$, where $R_\gamma\subset U_X$, approximates the set $I$ if
$
\lim_{\gamma \to \infty}\xi(I,R_\gamma) =0.
$
\medskip

(A.1)  {\it Prolongation}.\/ 
Let $z=(z_1, \dots, z_n)$ be a local coordinate system in a small neighborhood $\Cal V$ with origin at a {\it  fixed} point $\bold a\in \Cal V\subset U_X$. Let $\tilde \Phi_\bold a (z(p), \overline{z(p)})$ be the diastasic potential at $\bold a$. 
Let  $\bold b\in U_X$ be an arbitrary point. Let  
$$
\qquad \qquad\qquad \qquad I: u=u(s) \qquad (0\leq s\leq 1,\; u(0)= \bold a, u(1)= \bold b)
$$
 be a path joining $\bold a$ and $ \bold b$. 
\smallskip

(A.1.1)  {\it Prolongation along the path} $I$.\/ We say $\tilde \Phi_\bold a (z(p), \overline{z(p)})$ has a prolongation along $I$  if the following two conditions are satisfied:

\roster
\item"{(i)}" To every $s\in [0,1]$ there corresponds a functional element of diastasic potential $\tilde \Phi_{u(s)} (z(p), \overline{z(p)})$ at $u(s)$ ($u(s)$ is also called the center).

\item"{(ii)}" For every $s_0\in [0,1]$, we can take a suitable subarc $u:=u(s) \; (|s-s_0|\leq \epsilon, \epsilon >0)$ of $I$ contained in the domain of convergence of  $\tilde \Phi_{u(s_0)} (z(p), \overline{z(p)})$ such that every functional element $\tilde \Phi_{u(s)} (z(p), \overline{z(p)})$ with $|s-s_0|\leq \epsilon$ is a direct prolongation of $\tilde \Phi_{u(s_0)} (z(p), \overline{z(p)})$.
\endroster

The {\it direct prolongation} means the following. Suppose $\tilde \Phi_{\alpha_1} (z(p), \overline{z(p)})$ is defined on $\Cal V_1$ and $\tilde \Phi_{\alpha_2} (z(p), \overline{z(p)})$ is defined on $\Cal V_2$ ($\Cal V_1 \cap \Cal V_2 \not= \emptyset$). Then $\tilde \Phi_{\alpha_2} (z(p), \overline{z(p)})$ is the direct prolongation of $\tilde \Phi_{\alpha_1} (z(p), \overline{z(p)})$ if they coincide on $\Cal V_1 \cap \Cal V_2$.
Recall that the complexification allows us to consider the corresponding holomorphic function in place of the diastasic potential. It follows  a prolongation along $I$ is unique provided it exists. 

Let  $A = \{a_\nu\}$ ($\bold a, \bold b \in A$, $I=\bar A$) be a countable ordered dense subset of $I$. 
We would like to prolong $\tilde \Phi_\bold a (z(p), \overline{z(p)})$ along $I$ obtaining   the diastasic potential $\tilde \Phi_{a_\nu} (z(p), \overline{z(p)})$ of $\Lambda_\Cal L$ for each $a_\nu$.
We claim the prolongation along $I$ is possible. 
\smallskip

(A.1.2) Now, we will make use of the assumption $\pi_1(X)$ is large. We can assume $I$ is embedded in $X$ via $\tau$; otherwise, we could have replaced $X$ by $X_i$ for $i\gg0$. The set $A$ is a union of an increasing sequence of finite ordered subsets: 
$$
A_1 \subset A_2 \subset \cdots \subset A_\gamma \subset \cdots \subset A, \qquad \bold a, \bold b \in A_\gamma\; (\forall \gamma).
$$
We consider an arbitrary   $A_\gamma $ and the corresponding  set 
$$
\tau_i(A_\gamma)\subset \tau_i(I)\subset X_i \subset \bold P^{r_i},
$$ 
where $i$ is a sufficiently large integer and $r_i$ is an appropriate integer.  We apply Bertini's theorems to the linear system of curvilinear sections passing through $\tau_i(A_\gamma)$, i.e., the moving part of the system is a one-dimensional subscheme in $X_i$. 
We claim this linear system (and its inverse images) have no fixed components on $X_i$ for all $i\gg 0$.

Suppose, to the contrary, $W\subset X_i$ is a fixed  component. Then $W$ belongs to the linear span of $\tau_i(A_\gamma)\subset \bold P^{r_i}$. We move up along the tower (2.3.1.1).  For $j\gg i \gg 0$, the linear span of $\tau_j(A_\gamma)\subset \bold P^{r_j}$ will not contain $\tau^{-1}_{ji}(W)$. Hence the corresponding linear system on $X_j$ does not contain $\tau^{-1}_{ji}(W)$.  Therefore the linear system has no fixed components.

A priori, a general member of the system may have singularities at the base points  of  the system. However, we can always assume $\bold P^{r_i}$ is sufficiently large, and our system is sufficiently large as well.
Thus, the general member of the linear system on $X_j$ will be a connected nonsingular curve. Its inverse image on $U_X$ will be a connected open Riemann surface $R_\gamma$ by the Campana-Deligne theorem \cite{Kol, Theorem 2.14}.  These Riemann surfaces will approximate $I$ as $\gamma$ goes to infinity.  
\smallskip

(A.1.3)  
For every $\gamma$, the diastasic potential  of the induced metric on $R:= R_\gamma$ is the restriction of the corresponding diastasic potential of $U_X$, and  $\bold P_R:=\bold P_{R,\bold a}$ is a function on $R$.  

We replace the path $I$ and an arbitrary $A_\gamma\; (\gamma \gg 0)$ by a broken geodesic $\sigma_\gamma$ between the points $\bold a $ and $\bold b$. Namely, we replace the subpath of $I$ between two adjacent points of $A_\gamma$ by a geodesic on $U_X$. 
We also consider the corresponding broken geodesic $\rho_\gamma$ on $R_\gamma$.  Recall  (Section 2.2) that the diastasis approximates the square of the geodesic distance in the {\it small}. Hence $\rho_\gamma$ will be close to $\sigma_\gamma$ provided each $a_\nu$ is close to $a_{\nu+1}$,
and we get
$$
\lim_{\gamma \to \infty} \xi(I, \sigma_\gamma) =  \lim_{\gamma \to \infty} \xi(I, \rho_\gamma) = 0 \qquad (\sigma_\gamma \subset U_X, \; \rho_\gamma \subset R_\gamma).
$$
\smallskip

(A.1.4) Now, we will establish the prolongation along $I$. Assume we can prolong along $I\backslash \bold b$. We take a sufficiently small subarc $E\subset I$ in 
the domain  
$\Cal V_\bold b$ of 
$$
\tilde \Phi^{\Cal V_\bold b}_\bold b (z(p), \overline{z(p)}) := \tilde \Phi_\bold b (z(p), \overline{z(p)}).
$$
 Take a point   $\bold w\in E\backslash \bold b$ and its small neighborhood $\Cal V_\bold w \subset \Cal V_\bold b$ in $U_X$. We set  
$$
\tilde \Phi^{\Cal V_\bold w}_\bold b (z(p), \overline{z(p)}):= \tilde \Phi^{\Cal V_\bold b}_\bold b (z(p), \overline{z(p)})|\Cal V_\bold w,
$$
 more precisely,  $ \tilde \Phi^{\Cal V_\bold w}_\bold b (z(p), \overline{z(p)}) $ is a  real analytic function on $\Cal V_\bold w$ with center $\bold w$ (functional element)  obtained from the real analytic function $ \tilde \Phi^{\Cal V_\bold b}_\bold b (z(p), \overline{z(p)})$ on $\Cal V_\bold b$ with center $\bold b$.
For every $p\in \Cal V_\bold w$, we claim
$$
\Phi_{\bold w\bold b}(p) := \tilde \Phi^{\Cal V_\bold w}_\bold b (z(p), \overline{z(p)}) - \tilde \Phi^{\Cal V_\bold w}_\bold w (z(p), \overline{z(p)})= D_{U_X}^{\Cal V_\bold w}(\bold b, p) - D^{\Cal V_\bold w}_{U_X}(\bold w,p) =0, 
$$
where $D_{U_X}^{\Cal V_\bold w}(\bold b, p)$ is the real analytic function in $z(p), \overline{z(p)}$ on $\Cal V_\bold w$ with center $\bold w$ obtained from the real analytic function $D^{\Cal V_\bold b}_{U_X}(\bold b, p)$ in $z(p), \overline{z(p)}$ on $\Cal V_\bold b$ with center $\bold b$, and $D_{U_X}^{\Cal V_\bold w}(\bold w, p)$ is a real analytic function in $z(p), \overline{z(p)}$ on $\Cal V_\bold w$ with center $\bold w$. By the definition, $D_{U_X}^{\Cal V_\bold w}(\bold b, p)$ is a direct prolongation of $D^{\Cal V_\bold b}_{U_X}(\bold b, p)$. {\it A priori}\/, $D^{\Cal V_\bold w}_{U_X}(\bold w, p)$ is not a direct prolongation of $D^{\Cal V_\bold b}_{U_X}(\bold b, p)$.

Let $e\in \Cal V_\bold w$ be an arbitrary point.  We choose $\{R_\gamma\}$, as in (A.1.2), with an additional condition: $e, \bold w \in R_\gamma  \;(\forall \gamma)$. Then
$(\Phi_{\bold w\bold b}\!\big|{\!R_\gamma})(e) = 0$ for all $\gamma \gg 0$ because
$$
D^{\Cal V_\bold w}_{R_\gamma}(\bold b,e) - D^{\Cal V_\bold w}_{R_\gamma}(\bold w,e)=0 \qquad (\forall \gamma \gg 0)
$$
and the fundamental property of the diastasis (Section 2.2).
Here 
$$
D^{\Cal V_\bold w}_{R_\gamma}(\bold b,e)=D^{\Cal V_\bold w}_{U_X}(\bold b,e)\!\big|{\!R_\gamma} \quad {\text {and}}\quad D^{\Cal V_\bold w}_{R_\gamma}(\bold w,e)=D^{\Cal V_\bold w}_{U_X}(\bold w,e)\!\big|{\!R_\gamma} .
$$
It follows we can prolong $\tilde \Phi_\bold a (z(p), \bar z(p))$ along  $I$.
Since $U_X$ is simply connected, we obtain the desired function $\bold 
P_{U_X}\!:=\bold P_{U_X\!,\bold a}$ on $U_X$.

\proclaim{Lemma A (Prolongation)} Let $X\subset \bold P^r$ be a nonsingular connected
projective variety with large and residually finite $\pi_1(X)$. We assume the genus of
a general curvilinear section of $X$ is at least two and $U_X$ is equipped with the metric 
$\Lambda_\Cal L$ $($or, if $\pi_1(X)$ is nonamenable then $U_X$ is equipped with metric
$\beta$ as in $(5.1))$.  Then its diastasic potential\/ $\bold P_{U_X\!,\bold a}$ 
can be prolonged over $U_X$.

\endproclaim
\medskip

(A.2) In view of the Oka-Grauert-Narasimhan theorem (Grauert's version), it remains to verify that, for any real $\alpha$, the following set is relatively compact in $U_X$:
$$
E_\alpha :=\{u \in U_X \big |\; \bold P_U(u) < \alpha\}.
$$

Suppose $S\subset E_\alpha$ is an {\it infinite}\/ discrete subset without limit points in $U_X.$ Then we will derive a contradiction by showing that $\bold P_U$ is unbounded on $S$. Since $X$ is compact, $\tau(S)$ will be either a finite set or it will have a limit point. It suffices to replace $S$ by an {\it infinite}\/ set  $T_\alpha$ in the fiber of $\tau$ over a point $Q\in X$ and  
 show that $\bold P_U$ is unbounded on $T_\alpha$. If $\tau(S)$ has a limit point then $Q$ is such a point.

We consider a general curvilinear section $C\subset X$ through $Q$. Set
$
R_C:=\tau^{-1}(C) .
$ 
We obtain a connected open Riemann surface by the Campana-Deligne theorem \cite{Kol, Theorem 2.14}. 
          
By the fundamental property of the diastasis (Section 2.2), 
$\bold P_R = \bold P _U|_{R_C}$ where $\bold P_R$ is the corresponding diastasic  potentials on $R_C$. 
 One can find an infinite discrete subset $\tilde T_\alpha \subset \Cal F \subset \Delta$ (see (3.2)) whose image on $R_C$ will be close to the corresponding points of  $T_\alpha$. Moreover, $\tilde T_\alpha$ is approaching the boundary of $\Delta.$

We pick a point $\bold Q \in \Cal F$ such that its image on $X$ is close to $Q.$ We will identify $D_{\Lambda_R}(\bold Q, \cdot)$ with its inverse image on $\Cal F$.
We see that $D_{\Lambda_R}(\bold Q, \cdot)$ goes to infinity on $\tilde T_\alpha$
by considering the tower of coverings: 
$$
C\leftarrow \cdots \leftarrow C_i \leftarrow \cdots \leftarrow R_C,  
$$
 where $C_i \subset X_i$ (see (2.3.1.1)), and the diatasises of the corresponding  Bergman-type metrics of members of the tower restricted to the complement of hyperplane at infinity where they generate functions  as in (3.3.3). 

By Proposition 1, the diastasises increase as we move up in the tower. 
  So, $\bold P_R$ is unbounded on $\tilde T_\alpha$ and $T_\alpha,$ and $\bold P_U$ will be unbounded on $S$.

The contradiction proves the theorem.
\enddemo

\medskip

$(A.3)$ {\it  Remarks.}\/ Bogomolov and Katzarkov suggested   that the corresponding conjecture might fail in the case of nonresidually finite fundamental groups \cite{BK}.

A similar argument will establish the Shafarevich conjecture when $X$ is singular and $\pi_1(X)$ is residually finite and non-Abelian. Namely, let $X \hookrightarrow  \bold P^{r}$ be a  connected normal projective variety of dimension $n\geq 1$. Assume $\pi_1(X)$ is large, residually finite, and non-Abelian. Then its universal covering is a Stein space. The details will appear elsewhere.

\bigskip

{\it  Acknowledgment.} The author would like to thank 
J\'anos Koll\'ar and Gerard Washnitzer for helpful conversations, and Jean-Pierre Demailly who brought to my attention a theorem of M. Riesz. 
The author is grateful to Fedor Bogomolov, Fr\'ed\'eric Campana, Fabrizio Catanese, Phillip Griffiths, Curtis McMullen, Takeo Ohsawa, Peter Polyakov and Vyacheslav Shokurov  for their  emails.

\enddocument